\newcommand{\calA}{\mathcal{A}}
\newcommand{\calB}{\mathcal{B}}

\newcommand{\calE}{\mathcal{E}}
\newcommand{\calF}{\mathcal{F}}
\newcommand{\calG}{\mathcal{G}}
\newcommand{\calH}{\mathcal{H}}
\newcommand{\calK}{\mathcal{K}}
\newcommand{\calL}{\mathcal{L}}
\newcommand{\calM}{\mathcal{M}}

\newcommand{\calO}{\mathcal{O}}
\newcommand{\calP}{\mathcal{P}}
\newcommand{\calQ}{\mathcal{Q}}
\newcommand{\calR}{\mathcal{R}}
\newcommand{\calS}{\mathcal{S}}

\newcommand{\calJ}{\mathcal{J}}
\newcommand{\calI}{\mathcal{I}}

\newcommand{\calT}{\mathcal{T}}

\newcommand{\frakm}{\mathfrak{m}}
\newcommand{\frakc}{\mathfrak{c}}

\newcommand{\bbC}{\mathbb{C}}

\newcommand{\bbP}{\mathbb{P}}

\newcommand{\bbZ}{\mathbb{Z}}

\newcommand{\la}{\langle}
\newcommand{\ra}{\rangle}

\newcommand{\SL}{\textrm{SL}}

\newcommand{\Hom}{\text{Hom}}

\newcommand{\res}{\text{res}}
\newcommand{\st}{\text{s}}

\newcommand{\adj}{\text{adj}}

\newcommand{\Ar}{\text{Ar}}

\newcommand{\codim}{\text{codim}}
\newcommand{\ass}{\text{as}}
\newcommand{\se}{\text{ss}}

\newcommand{\Sym}{\text{Sym}}
\newcommand{\ome}{\Omega^1(\calA)}
\newcommand{\tome}{\tilde{\Omega}^1(\calA)}
\newcommand{\length}{\text{length}}

\documentclass[11pt]{amsart}
\usepackage{amscd, amsmath, amsthm, amssymb}

\newtheorem{theorem}{Theorem}[section]
\newtheorem{lemma}[theorem]{Lemma}
\newtheorem{proposition}[theorem]{Proposition}
\newtheorem{conjecture}[theorem]{Conjecture}
\newtheorem{corollary}[theorem]{Corollary}

\theoremstyle{definition}     
\newtheorem{definition}[theorem]{Definition}

\newtheorem{example}[theorem]{Example}

\theoremstyle{remark}
\newtheorem{remark}[theorem]{Remark}

\numberwithin{equation}{section}

\newcommand{\Quot}{\text{Quot}}

\newcommand{\Ext}{\text{Ext}}

\newcommand{\rank}{\text{rank}}

\newcommand{\PGL}{\textup{PGL}}
\newcommand{\op}{\calO_{\bbP^n}}

\newcommand{\Der}{\textup{Der}}

\begin{document}

\title[arrangements]
 {Logarithmic sheaves  attached to arrangements of hyperplanes}
\author[I. Dolgachev]{Igor V. Dolgachev}
\address{Department of Mathematics, University of Michigan, Ann Arbor, MI 48109, USA}
\email{idolga@umich.edu}
\thanks{Research  is partially supported by NSF grant
DMS-0245203}
\begin{abstract}
A reduced divisor on a nonsingular variety defines the sheaf of logarithmic 1-forms.  We introduce a certain coherent sheaf whose double dual coincides with this  sheaf.   It has some nice properties, for example, the residue exact sequence still holds even when the divisor is singular, and also it has a simple locally free resolution. We specialize to the case when the divisor is an arrangement of hyperplanes in projective space, and relate the properties of stability of the sheaf  with the combinatorics of the arrangement. We also extend a Torelli type theorem to non generic arrangements which allows one to reconstruct an arrangement from the sheaf attached to it. 
\end{abstract}
\maketitle



\section{Introduction} Any divisor $D$ on a nonsingular variety $X$ defines a sheaf of logarithmic differential forms $\Omega_X^1(\log D)$. Its equivalent definitions and many useful properties are discussed in a fundamental paper of K. Saito \cite{Sa}. This sheaf is locally free when $D$ is a strictly normal crossing divisor, and in this situation it is a part of the logarithmic De Rham complex used by P. Deligne to define the mixed Hodge structure on the cohomology of the complement $X\setminus D$. In the theory of hyperplane arrangements this sheaf arises when $D$ is a central arrangement of hyperplanes in $\bbC^{n+1}$. In exceptional situations this sheaf could be free (a free arrangement), for example, when the arrangement is a complex reflection arrangement. Many geometric properties of the vector bundle $\Omega_X^1(\log D)$ were studied in the case when $D$ is a generic arrangement of hyperplanes in $\bbP^n$ \cite{DK1}. Among these properties is a Torelli type theorem which asserts that two arrangements with isomorphic bundle of logarithmic 1-forms coincide unless they osculate a normal rational curve. In this paper we introduce and study a certain subsheaf $\tilde{\Omega}_X^1(\log D)$ of $\Omega_X^1(\log D)$. This sheaf contains as a subsheaf (and coincides with it in the case when the divisor $D$ is the union of normal irreducible divisors)  the sheaf of logarithmic differentials considered earlier in \cite{CHKS}. Its double dual is isomorphic to  $\Omega_X^1(\log D)$. Although $\Omega_X^1(\log D)$ could be locally free for very singular arrangements, e.g. when $n=2$ or for free arrangements, the sheaf $\tilde{\Omega}_X^1(\log D)$ is never locally free unless the divisor $D$ is locally formally isomorphic to a strictly normal crossing divisor. This disadvantage is compensated by some  good properties of this sheaf which $\Omega_X^1(\log D)$ does not posses in general. For example, one has always a residue exact sequence
$$0\to \Omega_X^1\to \tilde{\Omega}_X^1(\log D) \to \nu_*\calO_{D'}\to 0,$$
where $\nu:D'\to D$ is a resolution of singularities of $D$. Also, in the case when $D$ is an arrangement of $m$ hyperplanes in $\bbP^n$, the sheaf 
$\tilde{\Omega}_{\bbP^n}^1(\log D)$ admits a simple projective resolution  
$$0\to \calO_{\bbP^n}(-1)^{m-n-1} \to \calO_{\bbP^n}^{m-1} \to  \tilde{\Omega}_{\bbP^n}^1(\log D) \to 0.$$
In particular, its Chern polynomial does not depend on the combinatorics of the arrangement. This allows us to introduce the notion of a stable (resp. semi-stable, unstable) arrangement and define a map from the space of semi-stable arrangements to the moduli space of coherent torsion-free sheaves on $\bbP^n$ with fixed Chern numbers. All generic arrangements are semi-stable (and stable when $m \ge n+2$, and the Torelli Theorem mentioned above shows that the variety of semi-stable arrangements admits a birational morphism onto a subvariety of the moduli space of sheaves. We extend the Torelli theorem proving the injectivity on the set of semi-stable arrangements which contain a generic arrangement not osculating a normal rational curve and  conjecture that the same is true for all semi-stable arrangements whose dual configurations of points in $\check{\bbP}^n$ does not lie on the set of nonsingular points of a stable normal rational curve. We check the conjecture in the case of $\le 6$ lines in the plane. 

I am grateful to  Fabrizio Catanese,  Rob Lazarsfeld,  Mircea Musta\c{t}\u{a}, Giorgio Ottaviani  and Sergey Yuzvinsky for valuable remarks.

\section{Logarithmic 1-forms}
Let $X$ be a nonsingular $n$-dimensional algebraic variety over a field $k$ of characteristic 0 and $D$ be an effective Cartier divisor on $X$. Let $\Theta_{X/k}$ be the tangent sheaf on $X$ defined by $\Theta_{X/k}(U) = \Der_k(O_X(U))$, the $\calO_X(U)$-module of $k$-derivation of the coordinate ring $\calO_X(U)$.  Let $\phi_U = 0$ be a local equation of $D$ on $U$. Define a submodule 
$$\Theta_{X/k}(\log \phi_U) = \{\partial \in \Der_k(O_X(U)):\partial(\phi_U) \in (\phi_U)\}.$$
Since $\partial(a\phi_U) = \partial(a)\phi_U+a\partial(\phi_U) $, this definition does not depend on a choice of a local equation. Also, we have $\phi_U = g_{UV}\phi_V$ in $U\cap V$ and
$\partial(\phi_U) = \partial(g_{UV})\phi_V+g_{UV}\partial(\phi_V) $ shows that $\Theta_{X/k}(U)$ can be glued together to define a subsheaf $\Theta_{X/k}(\log D)$ of $\Theta_{X/k}$ and an exact sequence
\begin{equation}\label{eq}
0 \to \Theta_{X/k}(\log D)\to \Theta_{X/k} \to \calJ_D(D) \to 0,\end{equation}
where $\calJ_D$ is an ideal sheaf on $\calO_D$ generated in each $\calO_D(U)$ by $\partial(\phi_U), \partial\in  \Der_k(O_X(U)).$  In other words, 
$$\calJ_D = \text{Jacobian}(D)\cdot \calO_D,$$
where $\text{Jacobian}(D)$ is the {\it Jacobian ideal sheaf} in $\calO_X$ generated in each $\calO_X(U)$ by $\phi_U$ and $\partial(\phi_U), \partial \in \Der_k(O_X(U))$ (see \cite{La}, p.181). We set
$$\Omega_{X/k}^1(\log D): = \Theta_{X/k}(\log D)^* = \calH om_{X}(\Theta_{X/k},\calO_X)$$
and call it the {\it sheaf of logarithmic 1-forms} of $D$.
Since $\Theta_{X/k}$ is locally 
 free, dualizing \eqref{eq}, we get an exact sequence

\begin{equation}\label{eqqq}
0\to \Omega_{X/k}^1 \to \Omega_{X/k}^1(\log D)\to \calE xt_X^1(\calJ_D(D),\calO_X)\to 0.
\end{equation} 

 It follows from \eqref{eq} that 
$\text{depth}~ \Theta_{X/k}(\log D)_x \ge 2$ for any closed point. Thus the sheaf $\Theta_{X/k}(\log D)$ is reflexive, hence
$$\Theta_{X/k}(\log D)^{**} \cong \Omega_{X/k}^1(\log D)^* \cong \Theta_{X/k}(\log D).$$

Let $D^s$ be the closed subscheme of $D$ defined by the sheaf of ideals $\calJ_D$ so that $\calO_{D^s} = \calO_D/\calJ_D$. It is supported on the singular locus   of $D$. 

Consider the exact sequence
$$0\to \calJ_D(D) \to \calO_D(D) \to \calO_{D^s}(D) \to 0.$$
Applying the functor $\calH om_X(?,\calO_X)$ we get an exact sequence
$$0\to  \calE xt_X^1(\calO_D(D),\calO_X) \to \calE xt_X^1(\calJ_D(D),\calO_X) \to \calE xt_X^2(\calO_{D^s}(D),\calO_X) \to 0.$$

Let $\omega_Z$ denote the dualizing sheaf of a projective  Cohen-Macaulay algebraic  variety $Z$, the canonical sheaf $\calO_Z(K_Z)$ if $Z$ is nonsingular. By the Duality Theory,
$$\calE xt_X^1(\calO_X,\omega_X) \cong \omega_D \cong \omega_X\otimes_{\calO_X}\calO_D(D).$$
Therefore, 
\begin{equation}\label{relat}
\calE xt_X^1(\calO_D,\calO_X) \cong \calE xt_X^1(\calO_D,\omega_X)\otimes_{\calO_X}\omega_X^{-1} \cong  \calO_D(D).
\end{equation}

This proves the following: 

\begin{proposition}\label{omit} The sheaf $\calE xt_X^1(\calJ_D(D),\calO_X)$ from the exact sequence \eqref{eqqq} fits in the following exact sequence
$$0\to  \calO_D \to \calE xt_X^1(\calJ_D(D),\calO_X) \to \calE xt_X^2(\calO_{D^s}(D),\calO_X) \to 0.$$
\end{proposition}

It is known (see \cite{HL}, Proposition 1.1.6) that, for any coherent sheaf $\calF$ on $X$ supported on a closed subset of codimension $c$, 
\begin{equation}\label{hl}
\calE xt_X^q(\calF,\calO_X) = 0, \quad q < c.
\end{equation}

\begin{corollary} Assume that $\codim_XD^s \ge 3$. Then 
$$\calE xt_X^1(\calJ_D(D),\calO_X) \cong \calO_D,$$
and we have an exact sequence
$$0\to \Omega_{X/k}^1 \to \Omega_{X/k}^1(\log D) \to \calO_D\to 0.$$
\end{corollary}
\bigskip
Now let us recall the definition of the {\it adjoint ideal sheaf} $\adj(D)$ of $D$ (see \cite{La}, p. 179). Let $\mu:X'\to X$ be a birational morphism such that the proper inverse transform $D'$ of $D$ is nonsingular (a {\it log resolution} of $D$). Write $\mu^*(D) = D'+F$ for some divisor $F$ on $X'$ supported on the exceptional locus of $\mu$. We have
$$\adj(D) = \mu_*(K_{X'/X}-F),$$
where $K_{X'/X} = K_{X'}-\mu^*(K_X)$ is the relative canonical divisor of $\mu$.

\begin{lemma}\label{la} Let 
$$\frakc_D = \adj(D)\cdot \calO_D.$$
Then
\begin{itemize}
\item[(i)]  
$\calJ_D \subset \frakc_D$;
\item[(ii)] $\frakc_D\otimes \omega_D = \nu_*\omega_{D'}$;
\item[(iii)] $\adj(D) = \calO_X$ if and only if $D$ is normal and has at most rational singularities;
\item[(iv)] if $\nu:D'\to D$ is the normalization morphism with smooth $D'$, then $\frakc_D$ is the {\it conductor ideal sheaf}, i.e. the annulator sheaf of $\nu_*\calO_{D'}/\calO_D$.
\end{itemize}
\end{lemma}

\begin{proof} See \cite{La}, pp.179-181.
\end{proof}

\begin{proposition} The sheaf $\calE xt_X^1(\calJ_D(D),\calO_X)$ from the exact sequence \eqref{eqqq} fits in the following exact sequence
$$0\to \nu_*\calO_{D'}\to \calE xt_X^1(\calJ_D(D),\calO_X)  \to \calE xt_X^2((\frakc_D/\calJ_D)(D),\calO_X) \to 0.$$
\end{proposition}

\begin{proof}
It follows from part (iii) of Lemma \ref{la} that   $\frakc_D$ restricts to $\calO_{D}$  on the nonsingular locus of $D$, and so the sheaf $\calJ_D$.  This implies that $\frakc_D/\calJ_D$ is supported on the closed subset of codimension $\ge 2$ in $X$. Hence 
$\calE xt_X^1(\frakc_D/\calJ_D,\calO_X) = 0$. Also  
$\calE xt_X^2(\calO_D,\calO_X) = 0$ because $\calO_D$ has a locally free resolution of length 1 and $\calE xt_X^3(\calO_D/\frakc_D,\calO_X) = 0$ by \eqref{hl}.  This gives $\calE xt_X^2(\frakc_D,\calO_X) = 0$ and an exact sequence 
\begin{equation}\label{use}
\small{0\to \calE xt_X^1(\frakc_D(D),\calO_X) \to \calE xt_X^1(\calJ_D(D),\calO_X) \to 
\calE xt_X^2((\frakc_D/\calJ_D)(D),\calO_X) \to 0.}
\end{equation}
By adjunction formula, $\omega_D = \omega_X\otimes_{\calO_X}\omega_D(D).$ Applying part (ii) of Lemma \ref{la}, we get
$$\frakc_D(D) = \nu_*\omega_{D'}\otimes \omega_X^{-1}.$$
Hence
$$\calE xt_X^1(\frakc_D(D),\calO_X) = \calE xt_X^1(\nu_*\omega_{D'},\omega_X).$$
Applying Grothendieck's  Duality Theoren for a projective morphism (see \cite{Ha}, Theorem 11.1) together with Grauert-Riemenschneider's vanishing theorem $R^q\nu_*\omega_{D'} = 0,  q > 0,$ we have an isomorphism 
$$\calE xt_X^1(\nu_*\omega_{D'},\omega_X)  \cong \nu_*\calO_{D'}$$
Now the assertion follows from exact sequence \eqref{use}.

\end{proof}
\begin{definition} We set  $\tilde{\Omega}_{X/k}^1(\log D)$ to be the kernel of the composition map of sheaves
 $$\Omega_{X/k}^1(\log D) \to \calE xt_X^1(\calJ_D(D),\calO_X)\to  
 \calE xt_X^2((\frakc_D/\calJ_D)(D),\calO_X).$$
\end{definition}
 
  By definition, we have an exact sequences
\begin{equation}\label{res}
0\to \Omega_{X/k}^1 \to \tilde{\Omega}_{X/k}^1(\log D)\overset{\res}{\longrightarrow} \nu_*\calO_{D'}\to 0,
\end{equation} 
 We call this sequence the {\it residue} exact sequence. The reason for this name will be explained in the following example.
  
Also we have an exact sequence
\begin{equation}\label{second}
0\to \Omega_{X/k}^1(\log D)\to \tilde{\Omega}_{X/k}^1(\log D) \to\calE xt_X^2((\frakc_D/\calJ_D)(D),\calO_X) \to 0.
\end{equation}
Since $\calE xt_X^2((\frakc_D/\calJ_D)(D),\calO_X)$ is supported at a closed subset of codimension $\ge 2$, we have
$$\tilde{\Omega}_{X/k}^1(\log D)^{**} \cong \Omega_{X/k}^1(\log D)^{**} = \Omega_{X/k}^1(\log D).$$

\begin{proposition}\label{prop1} The following assertions are equivalent.
\begin{itemize}
\item[(i)] $\tilde{\Omega}_{X/k}^1(\log D) \cong \Omega_{X/k}^1(\log D);$
\item[(ii)] $(\frakc_D/\calJ_D)_x = \{0\}$ for any point $x\in D$ with $\dim \calO_{D,x} = 1$.
\end{itemize}

\end{proposition}

\begin{proof} Suppose (ii) holds. Then the sheaf $\frakc_D/\calJ_D$ is supported on a closed subset of $D$ of codimension $\ge 2$. By \eqref{hl}, $\calE xt_X^2((\frakc_D/\calJ_D)(D),\calO_X)  = 0$, and exact sequence \eqref{second} implies (i). Conversely, if (i) holds we have \\$\calE xt_X^2((\frakc_D/\calJ_D)(D),\calO_X)  = 0$. Passing to stalks at point $x\in D$ of codimension 1, we use that $\textup{Ext}_A^2(M,A) = 0$ for a module $M$ over a regular local ring of dimension 2 supported on the closed point implies $M = 0$. This easily follows from the fact that $\textup{Ext}_A^2(A/\frakm,A) \ne 0$, where $A/\frakm$ is the residue field of $A$. So, (ii) and (i) are equivalent. 
\end{proof}

\begin{definition} A divisor $D$ on $X$ is called a {\it normal crossing divisor} at a point $x\in D$ if $\calO_{D,x}$ is formally (or \'etale) isomorphic to the quotient of $\calO_{X,x}$ by an ideal generated by $t_1\ldots t_k$, where $t_1,\ldots,t_k$ is a subset of the set of local parameters in $\calO_{X,x}$. We say that $D$ is a normal crossing divisor in codimension $\ge k$ if $D$ is a normal crossing divisor at any point $x$ with $\dim \calO_{X,x} \le k$. A normal crossing divisor is a divisor which is normal crossing at each  point.
\end{definition}

 It is clear from the definition that a normal crossing divisor in codimension $\le 1$ is just a reduced divisor. A normal crossing divisor in codimension $\le 2$ is a divisor which is, in codimension $\le 2$,  formally isomorphic to the product of an affine space and an ordinary double point. 

\begin{corollary}\label{cor1} Suppose $D$ is a normal crossing in codimension $\le 2$. Then
$$\tilde{\Omega}_{X/k}^1(\log D) \cong \Omega_{X/k}^1(\log D).$$
The converse is true if for any  point  $x\in D$ of codimension 1 the local ring $\calO_{D,x}$ is locally (formally) can be defined in $\calO_{X,x}$ by an equation $u^a-v^b = 0$, where $u,v$ are local parameters of $\calO_{X,x}$. 
\end{corollary}

\begin{proof} If $D$ is a normal crossing in codimension $\le 2$ then 
a local computation shows that condition (ii) in Proposition \ref{prop1} is satisfied.  To prove the converse we may assume that $X$ is two-dimensional with local parameters $u,v$ at a point $x$ and $D$ is given by local equation $f(u,v) = u^a-v^b = 0$ at $x$. Then 
$$\length~\calO_{D,x}/\calJ_{D,x} = \length~\calO_{X,x}/(f_u',f_v',f) = \length~\calO_{X,x}/(u^{a-1},v^{b-1})  $$
$$= (a-1)(b-1).$$
Now we use a well-known  Jung-Milnor formula from the theory of curve singularities (see an algebraic proof in \cite{Ri})
\begin{equation}\label{milnor}\mu = 2\delta -r+1.
\end{equation}
Here 
$$\mu = \text{length}~\calO_{X,x}/(f_u',f_v'), \ \delta = \text{length}~\calO_{D,x}/\frakc_{D,x}$$
and $r$ is the number of local branches of $D$ at $x$. Write $a =md, b = nd$, where $(m,n) = 1$. Then 
$$u^a-v^b = (u^m)^d-(v^n)^d = \prod_{i=1}^d(u^m-\epsilon^iv^n),$$
where $\epsilon$ is a primitive $d$th root of unity. It follows that $d = r$ is the number of branches. By Proposition \ref{prop1}, $\delta = \mu$, hence by \eqref{milnor}, we get 
$$(a-1)(b-1) =  (md-1)(nd-1) = d-1.$$
This can happen only if $d=m=n = 1$ or $m =n = 1, d = 2$. In the first case $D$ is nonsingular at $x$. In the second case, $D$ is a normal crossing at $x$.
\end{proof}

\begin{remark} It follows from a result of Zariski \cite{Za} that the singularities $f = u^a-v^b = 0$ are characterized by the condition that $f\in (f_u',f_v')$, or equivalently, $\length~\calO_{D,x}/\calJ_{D,x} = \text{length}~\calO_{X,x}/(f_u',f_v')$
\end{remark}

\begin{definition} Let $Y$ be  a nonsingular subvariety of a nonsingular variety $X$ and $D$ be a reduced divisor on $X$. We say that $Y$  intersects $D$ transversally if $\calT or_1^{X}(\calO_Y,\calO_D) = 0$ and for any resolution of singularities $f:D'\to D$ the morphism 
$D'\times_XY\to D\times_XY = Y\cap D$ is a resolution of singularities. 
\end{definition}

\begin{proposition}\label{diagram} Let $Y$ be  a nonsingular subvariety $Y$ of $X$ with the sheaf of ideals $\calI$. Assume that $Y$ intersects transversally $D$. There is an exact sequence
$$0\to \calI/\calI^2 \to \Omega_{X/k}^1(\log D)\otimes_{\calO_X} \calO_Y \to \Omega_{Y/k}^1(\log D\cap Y) \to 0.$$
\end{proposition}

\begin{proof} We have a standard exact sequence
\begin{equation}\label{standard}
0\to \calI/\calI^2 \to \Omega_{X/k}^1\otimes_{\calO_X} \calO_Y \to \Omega_{Y/k}^1 \to 0.
\end{equation}
Consider the residue exact sequence for $(X,D)$ and tensor it with $\calO_Y$. Using the condition $\calT or_1^{X}(\calO_Y,\calO_D) = 0$, we get an exact sequence
$$0\to \Omega_{X/k}^1\otimes_{\calO_X} \calO_Y\to \Omega_{X/k}^1(\log D)\otimes_{\calO_X} \calO_Y
\to \nu_*\calO_{D'}\otimes_{\calO_X} \calO_Y \to 0.$$
Now consider the following commutative diagram

\[\small{\begin{array}{ccccccccccc}
&{}&0&&0&&0&&0&
\\
&&\Big\uparrow&{}&\Big\uparrow&{}&\Big\uparrow&{}&\Big\uparrow&\\
0&\longrightarrow&\calP&\longrightarrow&i^*(\nu_* \calO_{D'})&\longrightarrow& \pi_*\calO_{(D\cap Y)'}&\longrightarrow& \calQ&\longrightarrow&0\\
&&\Big\uparrow&{}&\Big\uparrow&{}&\Big\uparrow&{}&\Big\uparrow\\
0&\longrightarrow&\calR&\longrightarrow&i^*\Omega_{X/k}^1(\log D)&\longrightarrow& \Omega_{Y}^1(\log Y\cap D) &\longrightarrow& \calS&\longrightarrow& 0\\
&&\Big\uparrow&{}&\Big\uparrow&{}&\Big\uparrow&{}&\Big\uparrow\\
0&\longrightarrow&\calI/\calI^2&\longrightarrow&i^*\Omega_{X/k}^1&\longrightarrow&\Omega_{Y/k}^1 &\longrightarrow& 0&\longrightarrow& 0\\
&&\Big\uparrow&{}&\Big\uparrow&{}&\Big\uparrow&{}&\Big\uparrow\\
&{}&0&&0&&0&&0&
\end{array}}
\]

 Here $i:Y\hookrightarrow X$ is the inclusion morphism, and $\pi:(D\cap Y)'\to D\cap Y$ is a resolution of singularities which we can choose to be a composition of a resolution of singularities of $D'\times_XY$ and the projection $D'\times_XY\to D\times_XY = D\cap Y$. The middle horizontal exact sequence is obtained by dualizing a natural homomorphism
 $$ \Theta_{Y/k}(\log D\cap Y) \to \Theta_{X/k}(\log D)\otimes_{\calO_X}\calO_Y. $$
 In the row above it, we have a natural morphism of sheaves 
 $$\alpha:\nu_* \calO_{D'}\otimes_{\calO_{X}}\calO_{Y} \to   \nu_*\calO_{(D\cap Y)'}$$ which is the composition of an isomorphism $\nu_* \calO_{D'}\otimes_{\calO_{X}}\calO_{Y} \to   \nu_*\calO_{D'\times_XY}$ and a natural morphism  $\nu_*\calO_{(D'\times_XY)}\to \pi_*\calO_{(D\cap Y)'}$. By the transversality assumption, $D'\times_XY \cong (D\cap Y)'$, hence $\alpha$ is an isomorphism. This implies that $\calP = \calQ = 0$ and the assertion follows.
 \end{proof}

\begin{example}\label{normal} In the case when  $D$ is a {\it strictly normal crossing divisor}, i.e. the union of smooth divisors $D_i, i = 1,\ldots,m, $ which intersect transversally at each point, the sheaf $\Omega_{X/k}^1(\log D)$ and its exterior powers $\Omega_{X/k}^r(\log D)$ are well-known tools for defining the mixed Hodge structure on the complement $X\setminus D$. The sheaf $\Omega_{X/k}^1(\log D)$ is isomorphic to a subsheaf of the sheaf of rational differentials with poles on $D_i$ of order at most one.  If $z_i = 0, i = 1,\ldots,s,$  is a local equation of $D_i$ at  a point $x$ in the intersection $D_1\cap\ldots \cap D_s$, then $\Omega_{X/k}^1(\log D)$ is  locally free at $x$ and is generated in an open neighborhood of $x$ by meromorphic differential forms $d\log z_1,\ldots, d\log z_s, dz_{s+1},\ldots,dz_n$. Let $\epsilon_i:D_i\to X$ be the closed embedding. The map of sheaves 
$$\res:\tilde{\Omega}_{X/k}^1(\log D) \to \nu_*\calO_{D'} \cong \bigoplus_{i=1}^m\epsilon_{i*}\calO_{D_i}$$
is given by the residue map
$$\res(\sum_{i=1}^s a_id\log z_i+\sum_{s+1}^n b_idz_i) = (a_1+(z_1),\ldots, a_s +(z_s), 0,\ldots,0).$$
Since a normal crossing divisor is locally formally isomorphic to a simple normal crossing divisor, it follows that the sheaf $\Omega_{\bbP^n}^1(\log D)$ is locally free if $D$ is a normal crossing divisor.
\end{example}

\section{The logarithmic sheaves of a hyperplane arrangement}
This is a special case of the construction from the previous section. First we assume that $X$ is the projective space $\bbP^n$ over $k$ and $D$ is a hypersurface $V(f)$, where $f$ is a homogeneous element of degree $m$ in the polynomial algebra $S = k[T_0,\ldots,T_n]$.  Let 
$$\Omega_{S/k}^1 = SdT_0+\ldots+SdT_n \cong S(-1)^{n+1}$$ 
 and 
$$\Der_{S/k} = 
S\tfrac{\partial}{\partial T_0}+\ldots+S\tfrac{\partial}{\partial T_n} \cong S(1)^{n+1}$$
be the graded $S$-module of differentials and the graded $S$-module of derivations, dual to each other.  Recall that $S(a)_i = S_{a+i}$. Let
$E = \sum_{i=0}^n T_i\tfrac{\partial}{\partial T_i}$
be the Euler derivation. It defines a homomorphism of $E:\Omega_{S/k}^1\to S$ of graded modules. Let $\bar{\Omega}_{S/k}$ be its kernel. The corresponding sheaf on $\bbP^n$ is the sheaf $\Omega_{\bbP^n}^1$ of regular differential 1-forms. Its dual is the tangent sheaf 
$\Theta_{\bbP^n}$ associated to the cokernel of the homomorphism $S\to \Der_{S/k}, a\mapsto aE$. 
Let
$$\Der_{S/k}(\log f) = \{\partial\in \Der_{S/k}:\partial(f) \in (f)\}.$$
Obviously, $E\in  \Der_{S/k}(\log f)$. For any $\partial\in \Der_{S/k}$,  there exists a unique $p\in S$ such that $\alpha(\partial-pE) = 0$. Thus 
$$\Der_{S/k}(\log f) = SE\oplus \Der_{S/k}^0,$$
where $\Der_{S/k}^0$ is the kernel of the map $\Der_{S/k}\to S(m), \partial \mapsto \partial(f).$
Clearly, 
$$\widetilde{\Der}_{S/k}^0 \cong \Theta_{\bbP^n}(\log V(f)),$$
where $\tilde{}$ denotes the  sheaf associated to a graded $S$-module. Since $f\in J_f$, the ideal sheaf  $\tilde{J_f}$ on $\bbP^n$ can be considered as an ideal sheaf in $V(f)$ and it coincides with $\calJ_{V(f)}$ defined in the previous section. 

\bigskip
From now on we will consider the case when $f= f_1\cdots f_m$ is the product of distinct linear forms. The divisor $\calA =V(f)$ is called an {\it arrangement of hyperplanes}. We set
$$\Omega^1(\calA):= \Omega_{\bbP^n}^1(\log \calA),\quad \tilde{\Omega}^1(\calA):= \tilde{\Omega}_{\bbP^n}^1(\log \calA).$$
 It is customary in the theory of hyperplane arrangements  to grade $\Omega_{S/k}^1$ and its dual by assigning the grade zero to each $dT_i$ and $\frac{\partial}{\partial T_i}$. So their sheaf of logarithmic differentials is equal to $\Omega^1(\calA)(1).$

Let $L_i = V(f_i), i =1,\ldots,m,$ so that $V(f) = L_1\cup \ldots \cup L_m$. The normalization of $V(f)$ is isomorphic to the disjoint union of the $L_i$'s. We have 
\begin{equation}\label{for}
\nu_*\calO_{\calA'} = \bigoplus_{i=1}^m \epsilon_{i\ast}\calO_{L_i},
\end{equation}
where $\epsilon_i:L_i\hookrightarrow \bbP^n$ is the inclusion morphism. Since $\omega_{L_i} =  \calO_{L_i}(-n)$, we have
$$\nu_*\omega_{\calA'} =  \nu_*\nu^*\calO_{\bbP^n}(-n) = (\nu_*\calO_{\calA'})(-n)
= \bigoplus_{i=1}^m \epsilon_{i\ast}\calO_{L_i}(-n).$$
Thus 
$$\frakc_{\calA} = \nu_*\omega_{\calA'}\otimes \omega_{\bbP^n}^{-1} =
\bigl(\bigoplus_{i=1}^m \epsilon_{i\ast}\calO_{L_i}(-n)\bigl)\otimes \calO_{\bbP^n}(n+1) = \bigoplus_{i=1}^m \epsilon_{i\ast}\calO_{L_i}(1).$$
The following exact sequences are just the exact sequences \eqref{res} and \eqref{second} rewritten in our special situation
\begin{equation}\label{res'}
0 \to \Omega_{\bbP^n}^1 \to \tilde{\Omega}^1(\calA)\overset{\text{res}}{\longrightarrow} \bigoplus_{i=1}^m \epsilon_{i\ast}\calO_{L_i}\to 0,
\end{equation}
\begin{equation}\label{second'} 
0 \to \tilde{\Omega}^1(\calA) \to \Omega^1(\calA) \to \calE xt_{\bbP^n}^2((\frakc_{\calA}/\calJ_{\calA})(d+1),\op) \to 0.\end{equation}

\begin{theorem}\label{above} Assume $m\ge n+2$. The sheaf $\tilde{\Omega}^1(\calA)$ admits a projective resolution
$$0\to \op(-1)^{m-n-1} \to \op^m \to \tilde{\Omega}^1(\calA)\to 0.$$
\end{theorem}

\begin{proof} Let $i:\bbP^n \to \bbP^{m-1}$ be the closed embedding defined by 
$(t_0,\ldots,t_n) \mapsto (f_1,\ldots,f_m).$ Let $z_0,\ldots,z_{m-1}$ be projective coordinates in $\bbP^{m-1}$ and $\calB$ be the arrangement of the coordinate hyperplanes. Obviously, $i^*(\calB)  = \calA. $ We apply Proposition \ref{diagram}. 
 Formula \eqref{for} allows us to check the transversality condition.  Thus we have an exact sequence
 $$0\to \calI/\calI^2 \to i^*\Omega_{\bbP^{m-1}}^1(\log V(z)) \to \Omega^1(\calA)\to 0.$$
 The ideal sheaf $\calI$ of $i(\bbP^n)$ in $\bbP^{m-1}$ is associated to a free $k[z_0,\ldots,z_{m-1}]$-module generated by the subspace of linear polynomials spanned by $m-1-n$ linear independent linear relations between the functions $f_1,\ldots,f_m$. Thus 
 $$\calI/\calI^2 \cong \op^{m-n-1}(-1).$$ It is easy to check that 
 $$\Omega_{\bbP^{m-1}}^1(\calB) \cong \calO_{\bbP^{m-1}}^{m-1}$$
 (see \cite{DK1}, Proposition 2.10).  
 \end{proof}
 
Recall that an arrangement $\calA$ is called a {\it generic arrangement} if it is a simple normal crossing divisor.

\begin{proposition} The following assertions are equivalent
\begin{itemize}
\item[(i)] $\tilde{\Omega}^1(\calA)$ is locally free;
\item[(ii)] $\calA$ is a generic arrangement.
\end{itemize}
\end{proposition}

\begin{proof} It follows from Example \ref{normal} that (ii) implies (i). Assume (i) holds. Applying the residue exact sequence \eqref{res'}, we find that the sheaf $\nu_*\calO_{\calA'}$ is locally generated by $n$ elements. Suppose $\calA$ is not a normal crossing divisor. Then there exists a closed point $x\in \bbP^n$ such that there are $s > n$ hyperplanes $L_i$ passing through $x$.  Without loss of generality we may assume that $x = (1,0,\ldots,0)$ and the hyperplanes are given by linear equations $g_1,\ldots,g_{m}$ in inhomogeneous coordinates $z_1,\ldots,z_n$.  By \eqref{for}
$$(\nu_*\calO_{\calA'})_x \cong \bigoplus_{i=1}^s (k[z_1,\ldots,z_s]/(g_i))_{(z_1,\ldots,z_n)}.$$
We have a surjection $\calO_{X,x}^n \to (\nu_*\calO_{\calA'})_x$. After tensoring with \\$k[z_1,\ldots,z_n]_{(z_1,\ldots,z_n)}/(z_1,\ldots,z_n)$, we get a surjection of vector spaces $k^n \to k^s$. This contradiction proves the assertion.
\end{proof} 

\begin{proposition} The following assertions are equivalent
\begin{itemize}
\item[(i)] $\tilde{\Omega}^1(\calA) \cong \Omega^1(\calA)$;
\item[(ii)] $\calA$ is a normal crossing  divisor in codimension $\le 2$ (in $X$).
\end{itemize}
\end{proposition}

\begin{proof} This follows from Corollary \ref{cor1} since, locally in codimension 2,  the divisor $D$ can be written by equation $u^a-v^a = 0$, where $a$ is the number of hyperplanes in the arrangement $\calA$ intersecting along a codimension 2 subspace. 
\end{proof}

\begin{corollary}\label{known} Suppose $\calA$ is a normal crossing divisor in codimension $\le 2$. The following properties are equivalent
\begin{itemize}
\item[(i)] $\Omega^1(\calA)$ is locally free;
\item[(ii)]   $\calA$ is a generic arrangement.
\end{itemize}
\end{corollary}

\begin{remark} Recall that an arrangement $\calA$ is called {\it free} if the $S$-module $\Der_{S/k}(\log V(f))$ is free. Also $V(f)$ is called {\it locally free} if  the sheaf $\Omega^1(\calA)$ is locally free. Of course, a free divisor is locally free but the converse is not true in general.  If $n= 1$ any divisor is free but already in dimension 2  any reduced  divisor is locally free but not necessary free.  The assertion from Corollary \ref{known} follows from \cite{Zi} or \cite{Yu}, where it is proven that a free arrangement which is normal crossing in codimension $\le 2$ is a Boolean arrangement (i.e. consists of $n+1$ linear independent hyperplanes). For any $X$ from the lattice of the arrangement one considers the arrangement $\calA_X$ of hyperplanes which contain $X$. It is known that an arrangement is locally free if and only if each $\calA_X$ is free. The arrangement $\calA$ is normal crossing  if and only if each $\calA_X$ is Boolean. Another simple proof of this fact follows easily from \cite{MS}, where the Chern polynomial of $\ome$ is computed for a locally free arrangement (see  \eqref{chern}).
\end{remark}

\section{Stability of Steiner sheaves}

A coherent torsion-free sheaf $\calF$ on $\bbP^n$ with a projective resolution
$$0\to \op(d)^a \to \op(d+1)^b \to \calF \to 0, \quad 0<a<b,$$
 is called   a {\it Steiner sheaf} (see \cite{DK1}). 
 
Assume $m\ge n+2$. It follows from Theorem \ref{above} that the sheaf $\calF = \tilde{\Omega}^1(\calA)$ is a Steiner sheaf with the projective resolution
\begin{equation}\label{stres}
0\to \op(-1)^{m-n-1} \to \op^{m-1} \to \calF \to 0.
\end{equation}
Let $\bbP^n = \bbP(V)$ for some vector space $V$
$$U = H^0(\bbP(V), \calF\otimes \Omega_{\bbP(V)}^1(1)), \quad W = H^0(\bbP(V),\calF).$$
One identifies $U$ with $H^0(\bbP^n,\op(-1)^{m-n-1})$ by tensoring \eqref{stres} with $\Omega_{\bbP(V)}^1(1)$ and using the natural isomorphism $H^1(\bbP^n,\Omega_{\bbP^n}^1) \cong k$. Also one identifies $W$ with $H^0(\bbP^n,\op^{m-1})$. The map of sheaves $\op(-1)^{m-n-1} \to \op^{m-1}$ is defined by an injective linear map
$$t: V\to \Hom(U, W).$$
Conversely, one can reconstruct $\calF$ from such a map as the differential $d_{-1,0}$ in the Beilinson spectral sequence (see \cite{OSS}). 

In our situation when $\calF = \tilde{\Omega}^1(\calA)$, the proof of Theorem \ref{above} shows that   $U$ is isomorphic the subspace of $k^{m}$ which consists of relations between $f_i$'s, $W$ is isomorphic to the subspace of $k^{m}$ equal to the kernel of the map $(a_1,\ldots,a_m) \to \sum a_i$. The linear map $t$ is defined by the formula
\begin{equation}\label{tensor}
t(v)((a_1,\ldots,a_m)) = \bigl(a_1f_1(v),\ldots,a_mf_m(v)\bigr).
\end{equation}
(cf. \cite{DK1}). We will refer to $t_\calA:= t$ as the defining tensor of $\tilde{\Omega}^1(\calA)$. It could be considered as an element of the space $U^*\otimes V^*\otimes W$ and hence defines a divisor of multi-degree $(1,1,1)$ on $\bbP(U)\times \bbP(V)\times \bbP(W^*)$.  We say that $t_\calA$ is non-degenerate, if the  divisor is a nonsingular subvariety. The following proposition follows  easily from the definition. 

\begin{proposition}  $\tilde{\Omega}^1(\calA)$ is locally free if and only if $t_\calA$ is a non-degenerate tensor.\end{proposition}

Let $\calF$ be a torsion-free sheaf on $\bbP^n$. We identify its Chern classes with integers. It follows from \eqref{stres} that the Steiner sheaf $\tilde{\Omega}^1(\calA)$ has the Chern polynomial 
\begin{equation}
c_t(\tilde{\Omega}^1(\calA)) = 1/(1-t)^{m-1-n} = (1+t+\ldots+t^n)^{m-1-n} \mod (t^{n+1}).\end{equation}
Twisting \eqref{stres} by $\op(1)$, we also get
\begin{equation}\label{my}
c_t(\tilde{\Omega}^1(\calA)(1))  = (1+t)^{m-1} \mod (t^{n+1}) = \sum_{i=0}^nc_i(\Omega^1(\calA)t^i(1+t)^{n-i},
\end{equation}
where the last equality uses a well-known  relationship between the Chern polynomial of a sheaf and its Serre's twist.
On the other hand, if $\Omega^1(\calA)$ is locally free,  its Chern classes can be derived from \cite{MS}, Corollary 4.3:
\begin{equation}\label{chern}
P_{\calA}(t) = (1+t)c_t(\tilde{\Omega}^1(\calA)(1)),
\end{equation}
where $P_{\calA}(t)$ is the Poincar\'e polynomial of the arrangement
$$P_{\calA}(t) = \sum_{x\in \calL}\mu(x)(-t)^{\rank(x)}.$$
Here $\calL$ is the {\it lattice of the arrangement}, i.e. the partial ordered, by inclusion, set  of non-empty subsets 
$$L_I = L_{i_1}\cap \ldots,L_{i_s}, \quad I = \{i_1,\ldots,i_s\},$$
 $\mu:\calL\to \bbZ$ is the Moebius function of $\calL$ defined by
$$\mu(L_\emptyset) = 1,\  \ \mu(L_I) = -\sum_{L_I\subset L_J} \mu(L_J), $$
and $\rank(L_I) = \codim L_I$.

For a generic arrangement, we have
$P_\calA(t) = (1+t)^{m}$ and formulas \eqref{my} and \eqref{chern} agree.

Note that the  Poincar\'e polynomial $\Pi_{\calA}(t)$ of the corresponding central arrangement of affine hyperplanes in $k^{n+1}$ is related to ours $P_\calA(t)$ by the formula
$$\Pi_{\calA}(t) = P_\calA(t)-P_\calA(-1)(-t)^{n+1}.$$
\bigskip
\begin{example} Assume $n = 2$. Let $\calP$ be the set of singular points of $\calA$ (i.e. elements of $\calL$ of rank 2). We have $\mu(x) = s(x)-1$, where $s(x)$ is the number of lines through the point $x$. Then 
$$P_{\calA}(t) = 1+mt+\sum_{x\in \calP}(s(x)-1)t^2.$$
Using \eqref{chern}, we get  

\begin{eqnarray}\label{cherns}
c_1(\Omega^1(\calA)) & = & m-3, \\\notag
c_2(\Omega^1(\calA)) & = & \sum_{x\in \calP}(s(x)-1)-2m+3.
\end{eqnarray}
It follows from \eqref{second'} that 
$$c_1(\tilde{\Omega}^1(\calA)) = c_1(\Omega^1(\calA))$$
and
\begin{equation}\label{t}
c_2(\Omega^1(\calA)/\tilde{\Omega}^1(\calA)) = 
c_2(\Omega^1(\calA))-c_2(\tilde{\Omega}^1(\calA)) = \sum_{x\in \calP}(s(x)-1)-
\binom{m}{2}.\end{equation}
The second Chern class of a sheaf $\calT$ concentrated at a finite set of points is equal to 
$-h^0(\calT)$. Also, applying  Proposition \ref{above}, we get
\begin{equation}
h^0(\tilde{\Omega}^1(\calA)) = m-1,\quad  h^1(\Omega^1(\calA)) = 0.
\end{equation}
 Now \eqref{second'} gives
 \begin{equation}
 h^0(\Omega^1(\calA)) = m-1-\sum_{x\in \calP}(s(x)-1)+\tbinom{m}{2},\quad h^1(\Omega^1(\calA)) = 0.
 \end{equation}
 \end{example}

The $\rank~\calF$ is the rank of the vector bundle obtained by restriction to some open subset of $\bbP^n$. Recall that $\calF$ is called {\it semi-stable} (resp. {\it stable}) if for any proper subscheaf $\calF' \subset \calF$, 
$$\frac{h_{\calF'}(t)}{\rank~\calF'}  <  \frac{h_\calF(t)}{\rank~\calF},\ \ (\textrm{resp.} \frac{h_{\calF'}(t)}{\rank~\calF'} =  \frac{h_\calF(t))}{\rank~\calF}),$$
where $h_\calF(t) =\chi(\bbP^n,\calF(t))$ is the Hilbert polynomial of $\calF(t)$ and the inequality means the inequality between the values of the polynomials for $t >> 0$.

Comparing the coefficients at $t^{n-1}$, we see that  stability (resp. semi-stability) implies  {\it slope-stability}
$\mu(\calF')  < \mu(\calF)$ (resp.$\mu(\calF')  < \mu(\calF))$, where $\mu(\calF) = \frac{c_1(\calF)}{\rank~\calF}$ is the {\it slope} of $\calF$. The slope-stability implies stability but slope-semi-stability does not imply semi-stability.
 
 In the case $n = 2$ and $\calF$ is of rank $r$ with Chern classes $c_1$ and $c_2$, we have 
 $$\frac{h_\calF(t)}{r} = \frac{1}{2}t^2+(\mu(\calF)+3)t +\frac{3}{2}\mu(\calF)+\frac{1}{2r}(c_1^2-2c_2)+1.$$
 This shows that $\mu(\calF) =\mu(\calF')$ implies stability (resp. semi-stability) only if 
 $\Delta(\calF) > \Delta(\calF')$ (resp. $=$), where 
 $$\Delta(\calF) = \frac{1}{r}(c_2-\frac{r-1}{2r}c_1^2) = \frac{1}{2r}(c_1^2-2c_2)+\frac{1}{2}\mu(\calF)^2$$
is the {\it discriminant} of $\calF$. 

It is known that there is a coarse  moduli space $\calM_{\bbP^n}(r;c_t)$ of torsion-free semi-stable sheaves of rank $r$ on $\bbP^n$ with fixed Chern polynomial $c_t$ (\cite{Ma}). It is a projective variety. If $n = r =2$, we have
\begin{equation}\label{dim}
\dim \calM_{\bbP^2}(2;c_1,c_2) =  4c_2-c_1^2-3,
\end{equation} 
if the open subset of the moduli space representing stable sheaves is not empty. If any semi-stable sheaf is stable (e.g., if $(c_1,r) = 1)$), 
then  $\calM_{\bbP^n}(r;c_t)$ is a fine moduli space.

\begin{proposition} Assume $n > 1$. Any Steiner vector bundle $\calE$  on $\bbP^n$ defined by an exact sequence
$$0\to \op(-1)^{m-n-1} \to \op^{m-1} \to \calE \to 0$$
is a stable bundle of rank $n$ with the Chern polynomial $c_t = (1-t)^{n-m+1}$.
\end{proposition}

\begin{proof} It is enough to show that $\calE$ is slope-stable. This was proven in \cite{BS}. 
\end{proof}

It follows from \cite{DK1}, Corollary 3.3, that Steiner bundles (twisted by $\op(1)$) form an open subset $\calS_{n,m}$ in an irreducible component of the moduli space $\calM_{\bbP^n}(n;(1+t)^{m-1})$. If $n = 2$,  
$$\dim \calS_{2,m} = m(m-4).$$
 The logarithmic bundles $\Omega^1(\calA)$ of generic arrangements on $\bbP^2$ depend on $nm$ parameters. One proves that the map from the variety of general arrangements of $m$ hyperplanes  to the moduli space of vector bundles on $\bbP^n$ is a birational morphism for $m\ge n+2$. This was proved first in \cite{DK1} for $m\ge 2n+3$ and improved later in \cite{Va}. Thus for $n =2$, only in the case $m= 6$ we get the equality of the dimensions.

Now let us consider the problem of stability of Steiner sheaves $\calF$ on $\bbP^n = \bbP(V)$, not necessary locally free. We assume that 
$$\rank~\calF = n,$$
hence $\calF$ is given by an exact sequence
$$0\to \op(-1)\otimes U \to \op\otimes W \to \calF \to 0,$$
where 
$U \cong H^0(\bbP^n,\calF\otimes \Omega_{\bbP^n}(1)), \  W \cong H^0(\bbP^n,\calF)$ and the sheaf $\calF$ is determined by a tensor $t:V\to \Hom(U,W)$. We fix vector spaces $U$ and $W$ of dimensions $m-1-n$ and $m-1$, respectively and consider  the triples $(\calF,a,b)$, where $\calF$ is a Steiner sheaf and
$a,b$ are isomorphisms from above. Each such triple (a {\it Steiner triple}) is represented by a tensor $t$ defining a point in $\bbP(U^*\otimes V^*\otimes W)$. The condition of non-degeneracy is defined by a non-vanishing of the hyperdeterminant. Recall from \cite{GKZ} that the dual variety of   $\bbP_k^{n_1}\otimes\ldots \otimes \bbP_k^{n_s}$, embedded by Segre, is a hypersurface if and only if $n_i \le \sum_{j\ne i}n_j$ for any $i$. A tensor $t\in V_1\otimes \ldots \otimes V_s$, where $\bbP^{n_i} = \bbP(V_i)$, defines a hyperplane section of the Segre variety. So, it is singular if only if the hyperderterminant (which is an element of $\otimes_{i=1}^sV_i^*)$ vanishes at $t$. In our case $n_1+1 =\dim U = m-1-n, n_2+1 = \dim V = n+1, n_3+1 = \dim W = m-1$, so  $n_1 = n_2+n_3-2n, n_2 = n_1+n_3+2(m-n-2), n_3 = n_1+n_2$. Thus the hyperdeterminant exists if $m\ge n+2$. 

Let 
$$X_{m,n} = \bbP(U^*\otimes V^*\otimes W)/\!/\SL(U)\times\SL(W).$$
We can also view $X_{m,n}$ as the GIT-quotient of the Grassmannian of $m-1-n$-subspaces in $V^*\otimes W$:
$$X_{m,n} = G(m-1-n,V^*\otimes W)/\!/\SL(W).$$

The following result describes the set of semi-stable points in $G(m-1-n,V^*\otimes W)$ with respect to the action of $\SL(W)$ (\cite{Ka}, \cite{Ca}).

\begin{proposition}\label{unstable} A subspace $E\in G(m-1-n,V^*\otimes W)$ is semi-stable (stable) if and only if for each proper linear subspace $W'\subset W$ we have
$$\frac{\dim E\cap (W'\otimes V^*)}{\dim W'} \le \frac{\dim L}{\dim W} \ (\text{resp.}\ <)$$
\end{proposition}

\begin{corollary} Let $(\calF,a,b)$ be a Steiner triple with the defining tensor $t\in U^*\otimes V^*\otimes W$. Assume that $\calF$ is slope semi-stable (resp. slope stable). Then the tensor $t$, considered as a point in $ G(m-1-n,V^*\otimes W)$ is stable (resp. semi-stable).
\end{corollary}

\begin{proof} Let $E \subset V^*\cap W$ considered as the image of $U$ under the map $t:U\to V^*\otimes W$ defined by $t$. Let $U' = t^{-1}(E\cap W') \subset U$. It gives an exact sequence
of sheaves
$$0\to \op\otimes U' \to \op(1)\otimes W' \to \calF' \to 0.$$
It is clear that $\calF'(-1)$ is a subsheaf of the Steiner sheaf $\calF$ with 
$$\mu(\calF'(-1)) = \frac{\dim U'}{\dim W'-\dim U'}.$$
Since $\calF$ is slope stable (resp. slope semi-stable), we have
$$\frac{\dim U'}{\dim W'-\dim U'} \le \mu(\calF) = \frac{\dim U}{\dim W-\dim U}, \ (\text{resp.}\ < ).$$
It is easy to see that this is equivalent to the condition of semi-stability (stability) from the previous proposition.
\end{proof}

\begin{remark}\label{converse} The validity of the converse of the assertion in the previous corollary is unknown. It is true in the case when $m = n+3$ and $n$ is odd (see \cite{Ca}).
\end{remark}

\begin{corollary}\label{cor} Let $\calA$ be an arrangement of $m$ hyperplanes in $\bbP^n$ and $\calL$ be its lattice. For any $x\in \calL$ let $s(x)$ denote the number of hyperplanes containing $x$ and let $r(x) =\rank (x)$. Assume that there exists $x\in \calL$ such that 
$$s(x) > \frac{m-1}{n}(r(x)-1)+1.$$
 Then the Steiner log-sheaf $\tilde{\Omega}^1(\calA)$ is unstable. If the equality holds, $\tilde{\Omega}^1(\calA)$ is not stable.
\end{corollary}

\begin{proof} Assume such $x = L_I$ with $r(x) = r$ exists. Without loss of generality we may assume that the hyperplanes containing $L_I$  are  the hyperplanes $L_i = V(f_i), i= 1,\ldots,s$ and $f_1,\ldots,f_{r}$ are linearly independent. This implies that, for any $i = r+1,\ldots,s$, we can write $f_i = \sum_{j=1}^{r}a_{ij}f_j$. The corresponding relations span a subspace $U'$ of $U$ of dimension $s-r$. By definition of the defining tensor of $\calA$, it maps $U'$ to the subspace $V^*\otimes W'$ of 
$V^*\otimes W \subset V^*\otimes k^{m}$ generated by 
$$(a_{r1}f_1,\ldots,a_{rr}f_{r},-f_{r+1},0,\ldots,0), \ldots, $$
$$(a_{s1}f_1,\ldots,a_{sr}f_{r},0,\ldots,0,-f_{s},0,\ldots,0).$$
 Thus, in the notation of Proposition \ref{unstable}, we have  $\dim W' = s-1$ and 
$\dim U' = s-r = \dim E\cap W\otimes V^*$ and
$$\frac{\dim E\cap (W'\otimes V^*)}{\dim W'} -  \frac{\dim E}{\dim W} $$
$$= \frac{s-r}{s-1}-\frac{m-1-n}{m-1} = \frac{sn-n-(m-1)(r-1))}{(m-1)(s-1)}.$$
By assumption, the last number is positive, hence $t$ is unstable. By Corollary \ref{cor}, the sheaf  $\tilde{\Omega}^1(\calA)$ is unstable.
\end{proof}

\begin{proposition}  The sheaf $\tilde{\Omega}^1(\calA)$ is slope stable (resp. slope semi-stable) if and only if the sheaf $\Omega^1(\calA)$ is slope-stable (resp. slope semi-stable).
\end{proposition}

\begin{proof} More generally, let
$$0\to \calF \to \calG \to \calK \to 0$$
be an exact sequence of sheaves with $\rank~\calK = 0.$ Since $c_1(\calK) = 0$ and $\rank~\calF = \rank~\calG$, we have
$$\mu(\calF) = \mu(\calG).$$
Let $\calF'$ be a subsheaf of $\calF$ with $\mu(\calF') > \mu(\calF)$, then $\calF'$ is a subsheaf of $\calG$ with $\mu(\calF') > \calG$. Thus $\calG$ is unstable if $\calF$ is. Conversely, if $\calG'$ is a subsheaf of $\calG$ with $\mu(\calG') > \mu(\calG)$, we take $\calF'$ to be the kernel of the projection to $\calK$. Since $c_1(\calK) = 0$, we have $\mu(\calF') = \mu(\calG') > \mu(G) = \mu(\calF)$. Hence $\calF$ is unstable if $\calG$ is. 
This shows that slope semi-stability of $\calF$ is equivalent to slope semi-stability of $\calG$.
A similar proof, with replacing strict inequalities with non strict inequalities proves that slope stability of $\calF$ is equivalent to slope stability of $\calG$. We apply this to our situation using exact sequence \eqref{second'}.
\end{proof}

\begin{definition} An arrangement of hyperplanes $\calA$ is called \emph{stable} (resp. \emph{semi-stable}, resp. \emph{unstable}) if the sheaf $\tilde{\Omega}^1(\calA)$, or, equivalently, the sheaf $\Omega^1(\calA)$ is stable (resp. semi-stable, resp. unstable).
\end{definition}

\begin{example}\label{ex2} Let $\calA$ be a free arrangement. In this case the module of differentials 
$\Omega_{S/k}^1(\log f)$ is free, hence isomorphic to a direct sum of modules of type $S(a_i)$.
This shows that 
\begin{equation}\label{freesheaf}
\Omega^1(\calA) \cong \bigoplus_{i=1}^n\op(a_i).
\end{equation}
Its slope is equal to $(a_1+\ldots+a_n)/n$. Let us assume that $a_1\le \ldots\le a_n$. Then
the inequality $a_n\ge (a_1+\ldots+a_n)/n$ shows that $\mu(\op(a_n))\ge \mu(\Omega^1(\calA) )$ with equality only in the case $a_1=\ldots = a_n$. Hence $\Omega^1(\calA)$ is  unstable unless $a_1=\ldots=a_n$ in which case it is semi-stable.
\end{example}

\begin{example}\label{ex3} Take $n = 2$. The only interesting $r$ is $r = 2$, i.e. $x$ is a point in $\bbP^2$. We get that 
$s(x) > \frac{m-1}{2}+1$ implies unstability. For example, if $m = 6$, we need 4 lines passing through $x$. One should compare it with an inductive sufficient condition for slope stability and slope semi-stability of the bundle $\Omega^1(\calA)$  from \cite{Sch}, Theorem 4.5.  Note that the condition  $s(x) \le 3$ for any $x$ with $\rank(x) = 2$ is not sufficient for semi-stability. The reflection arrangement of type $A_3$ (its dual set of points in $\check{\bbP}^2$ is the set of vertices of  a complete quadrilateral) is free. By \eqref{cherns},  
$c_t(\Omega^1(\calA))) =1+3t+2t^2 = (1+t)(1+2t)$, hence $a_1 = 1, a_2 = 2$ in \eqref{freesheaf}. This shows that $\Omega^1(\calA)$ is unstable. This also can be proved without appealing to the freeness of the arrangement. It is known (\cite{OSS}, p. 168) that a vector bundle $\calE$ on $\bbP^2$ is unstable if 
$$8\Delta(\calE) = 4c_2(\calE)-c_1(\calE) < 0.$$
By \eqref{cherns}, this is equivalent to the inequality
\begin{equation}\label{unst}
4 \sum_{x\in \calP}(s(x)-1)-(m-1)(m+3) < 0.
\end{equation}
In the case of $A_3$-arrangement, the left-hand-side is equal to $44-45< 0$, so the sheaf
$\Omega^1(\calA)$ is unstable.
\end{example}

\bigskip
Recall that for any arrangement $\calA$ in $\bbP^n = \bbP(V)$ there is the {\it associated arrangement} $\calA^{\ass}$ (defined only up to projective equivalence) in $\bbP^{m-n-2} = \bbP(U)$ (see \cite{DK1}). The corresponding sheaf $\tilde{\Omega}^1(\calA^{\ass})$ is the Steiner sheaf defined by the same tensor $t\in U^*\otimes V^*\otimes W$ with the role of $U$ and $V$ exchanged. 

For any arrangement one defines the subset $D(\calA)$ of the set of subsets of $\{1,\ldots,m\}$ of cardinality $n+1$ which consists of subsets $(i_0,\ldots,i_n)$ such that 
$V(f_{i_1})\cap \ldots \cap V(f_{i_{n+1}}) \ne \emptyset$. In terms of the matrix of coordinates of the functions $f_i$, this is just the set of vanishing minors of maximal order. It follows from \cite{DO}, Lemma 1, p. 37, that the map $I\mapsto \{1,\ldots,m\}\setminus I$ is a bijection between the sets $D(\calA)$ and $D(\calA^{\ass})$. In particular, $\calA$ is  generic if and only if $\calA^{\ass}$ is  generic. 

\begin{conjecture} $\tilde{\Omega}^1(\calA)$ is stable if and only if $\tilde{\Omega}^1(\calA)^{\ass}$ is stable.
\end{conjecture}

\section{Unstable hyperplanes}
 
Let $\Ar_{n,m}$ be the variety of arrangements of $m\ge n+2$ hyperplanes in $\bbP^n$. This is just an open Zariski subset of $(\check{\bbP}^n)^{m}/S_{m}$ or, equivalently, a locally closed subset of the projective space of polynomials of degree $m$ which consists of products of $m$ distinct linear polynomials. We denote by $\Ar_{n,m}^{\se}$ (resp. $\Ar_{n,m}^{\st}$) the subset of semi-stable (resp. stable) arrangements. Let $\calS_{n,m}$ be a connected component of the Maruyama moduli space $\calM_{\bbP^n}(n,(1-t)^{n-m+1})$ which contains Steiner vector bundles defined by exact sequence \eqref{stres}. We have a map
\begin{equation}\label{map}
\log:\Ar_{n,m}^{\se}\to \calS_{n,m}, \quad \calA \mapsto \tilde{\Omega}^1(\calA).
\end{equation}
It is known that  this map is injective on the subset of generic arrangements which do not osculate a normal rational  curve of degree $n$ (i.e. the corresponding points in the dual projective space do not lie on such a curve). This was first proven in \cite{DK1} in the case $m\ge 2n+3$ and reproved under a weaker assumption $m\ge n+2$ in \cite{Va}. The generic arrangements osculating a  normal rational curve are blown down to the locus of Schwarzenberger bundles. 

The main idea of Valles's proof is to reconstruct the hyperplanes from the arrangement as {\it unstable} hyperplanes of the sheaf $\tilde{\Omega}^1(\calA)$. 

\begin{definition} Let $\calF$ be a Steiner sheaf of rank $n$ on $\bbP^n$. A hyperplane $L$ is called an {\it unstable} hyperplane of $\calF$ if 
$$H^{n-1}(L,\calF(-n)|L) \ne \{0\}.$$
\end{definition}

Here we denote by $\calF|L$ the scheme-theoretical restriction, i.e. 
$$\calF|L = i^*\calF = F\otimes_{\op}\calO_L,$$
where $i:L\hookrightarrow \bbP^n$ is the inclusion map.

\begin{proposition} Let $L$ be a hyperplane from a hyperplane arrangement $\calA$. Then $L$ is an unstable hyperplane of the sheaf $\tilde{\Omega}^1(\calA)$.
\end{proposition}

\begin{proof} Without loss of generality we may assume that $L = L_1$. We use the residue exact sequence \eqref{res'}.  Tensoring it with $\calO_L$ we obtain  an exact sequence
\begin{equation}\label{exter}
0\to \calT or_1^{\bbP^n}(\calO_L,\calO_L) \overset{\alpha}{\to} \Omega_{\bbP^n}^1|L 
\to \tilde{\Omega}^1(\calA)|L\to (\epsilon_{1\ast}(\calO_L\oplus 
\bigoplus_{i=2}^{m}\calO_{L_t\cap L})  \to 0
\end{equation}
Consider the exact sequence 
$$0\to \op(-1) \to \op \to \calO_L \to 0$$
corresponding to the inclusion of the ideal sheaf of $L$ in $\op$. Tensoring it with $\calO_L$, we get an exact sequence
$$0\to \calT or_1^{\bbP^n}(\calO_L,\calO_L) \to \calO_L(-1) \to \calO_L\to \calO_L \to 0.$$
This shows that $ \calT or_1^{\bbP^n}(\calO_L,\calO_L) \cong \calO_L(-1)$. Using \eqref{standard}, it is easy to identify the cokernel of the map $\alpha$ with $\Omega_{L}^1$. Thus we get an exact sequence
$$0\to \Omega_L^1 \to \tilde{\Omega}^1(\calA)|L\to \epsilon_{1\ast}(\calO_L\oplus 
\bigoplus_{t=2}^{m}\calO_{L_t\cap L})  \to 0.$$
Twisting by $\calO_L(-n)$ and applying cohomology, we get a surjection
$$H^{n-1}(L,\tilde{\Omega}^1(\calA)(-n)|L) \to H^{n-1}(L,\calO_L(-n)\oplus 
\bigoplus_{i=2}^{m}\calO_{L_i\cap L}(-n)) =$$
$$ H^{n-1}(L,\calO_L(-n)) = k.$$
This proves the assertion. 
\end{proof}

\begin{lemma} Let $\calA'$ be the arrangement obtained from an arrangement $\calA$ of $m\ge n+3$ hyperplanes by deleting a hyperplane $L$. There exists an exact sequence
$$0\to \tilde{\Omega}^1(\calA') \to \tome \to \epsilon_*\calO_L \to 0.$$
\end{lemma}

\begin{proof} The assertion probably follows from the residue exact sequence without the assumption on $d$, but this requires the verification that $\res^{-1}(\epsilon_*\calO_L)$ is isomorphic to $\tilde{\Omega}^1(\calA')$, so we prefer to give a simpler proof. We use that $\tome$ and $\tilde{\Omega}^1(\calA')$ are Steiner sheaves. We have a commutative diagram

\[\begin{array}{ccccccccc}
&{}&0&&0&&0&&
\\
&&\big\uparrow&{}&\big\uparrow&{}&\big\uparrow&{}&\\
0&\longrightarrow&\op(-1)&\longrightarrow&\op&\longrightarrow& \calO_L &\longrightarrow& 0\\
&&\big\uparrow&{}&\big\uparrow&{}&\big\uparrow&{}&\\
0&\longrightarrow&\op(-1)^{m-n-1}&\longrightarrow&\op^{m-1}&\longrightarrow& \tome &\longrightarrow& 0\\
&&\big\uparrow&{}&\big\uparrow&{}&\big\uparrow&{}&\\
0&\longrightarrow&\op(-1)^{m-n-2}&\longrightarrow&\op^{m-2}&\longrightarrow& \tilde{\Omega}^1(\calA') &\longrightarrow& 0\\
&&\big\uparrow&{}&\big\uparrow&{}&\big\uparrow&{}&\\
&{}&0&&0&&0&&
\end{array}
\]
Here the top horizontal sequence is the exact sequence of the definition of the sheaf $\calO_L$. The first two vertical exact sequences are obtained from composing the defining tensor 
$t_\calA:V\to \Hom(U,W)$ of $\calA$ with the restriction map $\Hom(U,W)\to \Hom(U',W')$, where $t_{\calA'}:V\to \Hom(U',W')$ is the defining tensor of $\calA'$. The right vertical sequence is the needed exact sequence.
\end{proof}

\begin{proposition}\label{uno} Let $\calA'$ be the arrangement obtained from an arrangement $\calA$ by deleting a hyperplane $L$. Then $W(\tome) \subset W(\tilde{\Omega}^1(\calA'))\cup \{L\}.$
\end{proposition}

\begin{proof} It is enough to  that  any $L'\in W(\tome)\setminus \{L\}$ belongs to $W(\tilde{\Omega}^1(\calA'))$. Tensoring the exact sequence from the previous Lemma by $\calO_{L'}(-n)$ we get an exact sequence
$$0\to \tilde{\Omega}^1(\calA')(-n)|L' \to \tome(-n)|L' \to \epsilon_*\calO_{L'\cap L}(-n) \to 0.$$
Taking cohomology, we get a surjection
$$\tilde{\Omega}^1(\calA')(-n)|L' \to H^{n-1}(L',\tome|L').$$
This shows that $L'\in W(\tilde{\Omega}^1(\calA'))$ if $L'\in W(\tome).$
\end{proof}
 
 In the case of general arrangements this result is Proposition 2.1 from \cite{Va} and Theorem 3.13 from \cite{AO} (where the inclusion is taken in scheme-theoretical sense, see below).
 
 \begin{corollary}\label{new} Assume   $\calA = \calA'+L$, where $\calA'$ is an arrangement such that $W(\tilde{\Omega}^1(\calA'))$ consists of $m-1$ unstable hyperplanes. Then 
 $$W(\tome) = W(\tilde{\Omega}^1(\calA'))\cup \{L\}.$$
 \end{corollary}
 
 \begin{proof}  $W(\tilde{\Omega}^1(\calA'))$ consists of hyperplanes from $\calA'$. Thus 
 $W(\tilde{\Omega}^1(\calA'))\cup \{L\} \subset W(\tome)$.  By Proposition \ref{uno}, we have the opposite inclusion. 
 \end{proof}
 
The set $W(\calF)$ of unstable hyperpanes of a Steiner sheaf $\calF$ has a natural structure of a closed subscheme of the dual projective space $\check{\bbP}^n$ (see \cite{AO}). In fact, one can construct a closed subscheme of $\calS_{n,m}\subset \calS_{n,m}\times \check{\bbP}^n$ such that the projection  
$$p: \tilde{\calS}_{n,m} \to \calS_{n,m}$$
has fibres isomorphic to the varieties $W(\calF)$ under the projection to the second factor.
The image of $p_1$ is a proper closed subvariety. Let 
$$p': \widetilde{\Ar}_{n,m}^{\se} \to \Ar_{n,m}^{\se}$$
be the pull-back of the map $p$ with respect to the map $\log:\Ar_{n,m}^{\se}\to \calS_{n,m}$. We know that over an open subset of generic arrangements which do not osculate a normal rational curve, the map $q$ is an unramified cover of degree $m$. Over the locus  of generic arrangements osculating a normal rational curve the fibres are isomorphic to $\bbP_k^1$. It follows that there exists an open Zariski subset $U\subset Ar_{n,m}^{\se}$ containing generic arrangements not osculating a normal rational curve such that, for any $\calF\in U$, the scheme  $W(\calF)$ is a reduced $0$-dimensional and consists of $m$ points. 

\begin{definition} An arrangement $\calA$ of $m$ hyperplanes is called a \emph{Torelli arrangement} if 
$W(\tome)$ consists of $m$ hyperplanes of $\calA$. 
\end{definition}

\begin{theorem} Let $U$ be the subset of $\Ar_{n,m}^{\se}$ which consists of Torelli arrangements.  Then $U$ is an open susbet of $\Ar_{n,m}^{\se}$ and the map $\log:U\to \calS_{n,m}$ is injective.
\end{theorem}

 Examples of Torelli arrangements are generic arrangements of $m\ge n+2$ which does not osculate a normal rational curve in $\bbP^n$ \cite{Va}. It follows from By Proposition \ref{new} that any arrangement which contains a Torelli arrangement is a Torelli arrangement. In particular any arrangement which contains a generic  arrangement $\calA'$ with at least $n+2$ hyperplanes not osculating a normal rational curve is Torelli. 

\begin{conjecture}\label{two}
A semi-stable arrangement of $m\ge n+2$ hyperplanes in $\bbP^n$ is always Torelli unless 
the corresponding points in $\check{\bbP}^n$ lie on a stable normal rational curve of degree $n$.
\end{conjecture}

Recall that a stable normal rational curve in $\bbP^n$ is a connected reduced curve of arithmetic genus 0 and degree $n$ in $\bbP^n$. It is the union of smooth rational curves $C_1,\ldots,C_s$ of degrees 
$d_1,\ldots,d_s$ satisfying the following conditions
\begin{itemize}
\item [(i)] $n = d_1+\ldots+d_s$;
\item[(ii)]  each curve $C_i$ spans a subspace $\la C_i\ra  = \bbP(V_i)$ of $\bbP^n = \bbP(V)$ of dimension $d_i$;
\item[(iii)] $V = V_1+\ldots+V_s$.
\end{itemize} 
\bigskip

\section{Line arrangements}
Here we assume $n = 2$. Recall that a line $L$ is called a {\it jumping line} of a rank 2 vector bundle $\calE$ on $\bbP^2$  if the splitting type of the restriction of $\calE$ to $L$ is different from the splitting type of the restriction of $\calE$ to a general line in the plane. This means that 

$$\calE|L \not\cong\begin{cases}
      \calO_L(a)\oplus \calO_L(a)& \text{if} \ c_1(\calE) = 2a, \\
\calO_L(a)\oplus \calO_L(a-1)& \text{if}\  c_1(\calE) = 2a-1.
\end{cases}$$
Equivalently, this means that $H^1(\calE(-a-1)|L) \ne 0$  if $c_1(\calE) =2a$ and 
$H^1(\calE(-a)|L) \ne 0$  if $c_1(\calE) =2a-1$. 
It is easy to see that $H^1(\calE(-2)|L) = 0$ implies $H^1(\calE(-2-s)|L) = 0$ for any $s\ge 0$. In \cite{DK1} an unstable line of $\Omega^1(\calA)$ for a generic arrangement $\calA$ was called a {\it super-jumping line}. Note that the notions of an unstable line of $\ome$ and a jumping line of $\ome$ coincide only if $m = 5$ or $6$. The exact sequence \eqref{second'} shows that any unstable line of $\tilde{\Omega}^1(\calA)$ not passing through its singular locus is a  jumping line of  $\Omega^1(\calA)$. 

Let $\calM_{\bbP^2}(2;c_1,c_2)$ be the moduli space of semi-stable sheaves of rank 2 on $\bbP^2$ with fixed Chern classes $c_1,c_2$. If there exists a stable vector bundle  with these Chern classes (e.g. if $(c_1,c_2) = 1$) then it is an irreducible variety of dimension $4c_2-c_1^2-3$ (\cite{Ma},\cite{Ba},\cite{Hu}). Consider its boundary $\partial\calM_{\bbP^2}(2;c_1,c_2)$ formed by sheaves which are not locally free. 
For any sheaf $\calF$ from the boundary, the double dual sheaf $\calF^{**}$ is a semi-stable vector bundle with the same $c_1$ and $c_2(\calF^{**}) = c_2-\delta$ for some $\delta\ge 0$. Let $\calM_{\bbP^2}(2;c_1,c_2)^{\delta}$ be the subset of $\calM_{\bbP^2}(2;c_1,c_2)$ which parametrizes isomorphism classes of such sheaves (or, more precisely, the corresponding $S$-equivalence classes if the sheaves are not stable but semi-stable). Since all bundles with $c_1^2-4c_2 > 0$ are known to be unstable (see \cite{OSS}, p.168), 
$$\calM_{\bbP^2}(2;c_1,c_2)^{\delta} = \emptyset, \quad \delta > 4c_2-c_1^2.$$
Note that 
$$\partial\calM_{\bbP^2}(2;c_1,c_2) = \cup_{\delta > 0}\calM_{\bbP^2}(2;c_1,c_2)^{\delta}.$$
Let
$$0\to \calF \to \calF^{**} \to \calT \to 0,$$
be the canonical exact sequence corresponding to the natural inclusion $\calF \subset \calF^{**}.$ The sheaf $\calT$ is concentrated at the set of singular points of $\calF$. Let $\delta_x$ be the length of the $\calO_{\bbP^2,x}$-module $\calT_x$. Let
$$Z(\calF) = \sum_{x\in \bbP^2}\delta_xx \in \Sym^\delta(\bbP^2)$$
be the corresponding point of the symmetric product of the plane. The set-theoretical union
$$\calM_{\bbP^2}(2;c_1,c_2)^U = \coprod_{\delta\ge 0}\calM_{\bbP^2}(2;c_1,c_2-\delta)^{0}\times \Sym^\delta(\bbP^2)$$
has a structure of a projective algebraic variety and is called the Uhlenbeck compactification of the moduli space of semi-stable vector bundles $\calM_{\bbP^2}(c_1,c_2)^0$ (see \cite{Li}). The natural map
$$\calM_{\bbP^2}(2;c_1,c_2)\to \calM_{\bbP^2}(2;c_1,c_2)^U, \quad \calF\mapsto (\calF^{**},Z(\calF))$$
is a morphism of algebraic varieties. Its
fibre  over a point $Z = \sum \delta_xx$ is isomorphic to the product of punctual quotient schemes 
$\Quot(2\delta_x)$ parametrizing quotient sheaves of $\calO_{\bbP^2}^2$ concentrated at $x$ and of length $\delta_x$. It is an irreducible variety of  dimension $2\delta_x-1$. There is an open subset of $\calM_{\bbP^2}(2;c_1,c_2)^{U}$ corresponding to points $Z= \sum_x\delta_xx$ such that $\delta_x\le 1$. The pre-image of this set in $\calM_{\bbP^2}(2;c_1,c_2)^\delta$ is an open subset of dimension equal to $\dim \calM_{\bbP^2}(2;c_1,c_2-\delta)$. It projection to $\calM_{\bbP^2}(2;c_1,c_2-\delta)^0$ has fibres of dimension $3\delta$.

\bigskip
Now let us specialize to our situation. Consider the exact sequence \eqref{second'} 
$$0\to \tilde{\Omega}^1(\calA) \to \Omega^1(\calA) \to \calT \to 0,$$
where $\calT = \calE xt_{\bbP^2}^2(\frakc_{\calA}/\calJ_{\calA},\calO_{\bbP^2}).$ The stalks of $\frakc_\calA$ and $\calJ_\calA$ are easy to compute using the Jung-Milnor formula from the proof of Corollary \ref{cor1}. We have 
$$\length (\frakc_\calA/\calJ_\calA)_x = \tbinom{s(x)-1}{2}.$$
Since $ \calE xt_{\bbP^2}^2(k,\calO_{\bbP^2})\cong k$, this gives
\begin{equation}
\length~\calT_x  = \tbinom{s(x)-1}{2}. 
\end{equation}
We know from \eqref{cherns} that 
$$h^0(\calT) = \sum_{x\in \calP}\length~\calT_x = \tbinom{m}{2}-\sum_{x\in \calP}(s(x)-1).$$
This gives a well-known combinatorial formula
\begin{equation}
\tbinom{m}{2}-\sum_{x\in \calP}(s(x)-1) = \sum_{x\in \calP}\tbinom{s(x)-1}{2}.
\end{equation}
We set
$$\delta_x(\calA):= \tbinom{s(x)-1}{2},\quad \delta(\calA): = \sum_{x\in \calP}\delta_x(\calA).
$$
Note that $\delta(\calA) = 0$ if and only if $\calA$ is a generic arrangement. It follows from \eqref{cherns}, that the numbers $d$ and $\delta$  determine the Chern polynomial of $\Omega^1(\calA)$. Recall that the moduli space of Steiner sheaves $\calS_{2,m}$ is equal to the moduli space $\calM_{\bbP^2}(2;c_1,c_2)$, where $c_1 = m-3, c_2= \binom{m-2}{2}$. Let 
$\calS_{2,m}^\delta = \calM_{\bbP^2}(2;c_1,c_2)^\delta$. Let $\Ar_{2,m}^{\se}(\delta)$ be the set of semi-stable arrangements with fixed $\delta(\calA) = \delta$. The restriction of the map \eqref{map} to $\Ar_{2,m}^{\se}(\delta)$ defines a map
$$\log_\delta:\Ar_{2,m}^{\se}(\delta) \to \calS_{2,m}^\delta.$$
One can rewrite the condition of unstability from \eqref{unst} in the form
\begin{equation}
\Ar_{2,m}^{\se}(\delta) = \emptyset, \quad \delta > \frac{(m-3)(m-1)}{5}.
\end{equation}
We also know from above that
$$\codim_{\calS_{2,m}}(\calS_{2,m}^\delta) = \delta.$$
Also taking the double dual defines a morphism 
$$u_\delta:\calS_{2,m}^\delta \to \calM_{\bbP^2}(2;m-3,\tbinom{m-2}{2}-\delta),$$ 
The composition 
$$u_\delta\circ \log_\delta:\Ar_{2,m}^{\se}(\delta) \to  \calM_{\bbP^2}(2;m-3,\tbinom{m-2}{2}-\delta)$$
is just the map $\calA\mapsto \Omega^1(\calA)$. It is easy to see that $\Ar_{2,m}^{\se}(1)$ is  irreducible and of codimension 1 in $\Ar_{2,m}^{\se}$. However, $\Ar_{2,m}^{\se}(2)$ consists of two irreducible components, each  of codimension 2. I do not know neither the number of irreducible component of $\Ar_{2,m}^{\se}(\delta)$ not their codimension for arbitrary $d$ and $\delta$.

\begin{remark} It follows from Schenk's inductive criterion of semi-stability \cite{Sch} that 
all arrangements with $\delta(\calA) = 1$ are stable for $m \ge 6$.
\end{remark}
 
 \begin{example} Let $m = 4$. Here only generic arrangements are stable. The moduli space 
 $\calM_{\bbP^2}(2;1,1) \cong \calM_{\bbP^2}(2;-1,1)$ consists of one point, representing the sheaf $\Omega_{\bbP^2}^1(2).$ Thus
 $$\tome = \ome \cong \Omega_{\bbP^2}^1(2) \cong \Theta_{\bbP^2}(-1).$$
 The exact sequence
 $$0\to \calO_L(-1) \to \Omega_{\bbP^2}^1|L  \to \Omega_L \to 0$$
 shows that 
 $$H^1(L,\ome(-2)|L) \cong H^1(L,\Omega_{\bbP^2}^1|L) \cong H^1(L,\Omega_L^1)\cong  k.$$ Thus any line is an unstable line of $\ome$.
 \end{example}
 
 \begin{example} Let  $m = 5$. The moduli space $\calS_{2,5} = \calM_{\bbP^2}(2;2,3) \cong \calM_{\bbP^2}(2;0,2)$ is a 5-dimensional variety.   Its open subset $\calS_{2,5}^0$ representing vector bundles is isomorphic to an  open subset $U$ of $\bbP^5$. If we identify the letter with the space of curves of degree 2 in the dual plane, then $U$ is equal to the set of nonsingular conics and the isomorphism is defined by assigning to a vector bundle $\calE$ its set of jumping lines (see \cite{Ba}). The variety $\calM_{\bbP^2}(2;2,2) \cong \calM_{\bbP^2}(2;0,1)$ is 2-dimensional. A sheaf $\calF$ from $\calM_{\bbP^2}(2;2,2)$ is determined by an extension
$$0\to \calO_{\bbP^2}\to \calF \to \calI_A(2) \to 0,$$
where $\calI_A$ is the ideal sheaf of a 0-dimensional closed subscheme in the plane with $h^0(\calO_A) = 2$. It shows that  $h^0(\calF(-1)) \ne 0$, hence $\calF$ contains a subsheaf $\calO_{\bbP^2}(1)$ of slope 1. Since $\mu(\calF) = 1$, this shows that  $\calM_{\bbP^2}(2;2,2)$ represents the $S$-equivalence classes of semi-stable but not stable sheaves. Each such class consists of vector bundles represented uniquely (up to isomorphism) by an extension
\begin{equation}\label{ext5}
 0\to \calO_{\bbP^2}(1) \to \calE \to \calI_x(1) \to 0
\end{equation}
for some point $x$. The only  non-locally semi-stable sheaf  in this class is the sheaf $\calO_{\bbP^2}(1)\oplus \calI_x(1)$, where $x$ is a point.  

The variety   $\calM_{\bbP^2}(2;2,1) \cong \calM_{\bbP^2}(2;0,0)$ is a one point. It represents the $S$-equivalence class of the sheaf $\calO_{\bbP^2}(1)^2.$

Thus for a generic arrangement $\calA$ of 5 lines we have $\tilde{\Omega}^1(\calA) \cong 
\Omega^1(\calA)$ is the Schwarzenberger vector bundle with the curve of jumping lines equal to the unique nonsingular conic in the dual plane containing the five lines of the arrangement. The map $\Ar_{2,5}^{\se}(0) \to \calM_{\bbP^2}(2;2,3)^0 = U$ is a surjective map with 5-dimensional fibres. 

The set $\Ar_{2,5}^{\se}(1)$ consists of arrangement with one triple point. All these arrangements are semi-stable but not stable. The sheaf 
$\Omega^1(\calA)$ belongs to $\calM_{\bbP^2}(2;2,2)$ and is $S$-equivalent to the sheaf $\calO_{\bbP^2}(1)\oplus \calI_x(1)$, where $x$ is a point. Observe that the two lines, say $L_1,L_2$  of $\calA$ not passing through the triple point are jumping lines of $\tilde{\Omega}^1(\calA)$ and hence of $\Omega^1(\calA)$. The set of unstable lines of a sheaf given by an extension \eqref{ext5}  is equal to the set of lines passing through $x$. This shows that $x = L_1\cap L_2$.

Thus all arrangements with the same point of intersection of two lines $L_0$ and $L_2$ not passing through the triple point have  bundle $\Omega^1(\calA)$ given by extension \eqref{ext5}, where $x = L_0\cap L_1$. The sheaf 
$\tilde{\Omega}^1(\calA)$ determines $\Omega^1(\calA)$ as its double dual, and determines the triple point $y$, as its singular point. So it determines a reducible conic in the dual plane, union of  the line dual to the triple point and the line dual to the  point $L_0\cap L_1$. All arrangements defining the same conic have the same $S$-equivalence class of the sheaf $\tilde{\Omega}^1(\calA)$. It is represented by the sheaf $\calI_x(1)\oplus \calI_y(1)$. Since $\Ext_{\bbP^2}^1(\calI_x(1),\calI_y(1)) \cong k$ if $x\ne y$, we obtain that there is a unique nontrivial extension class of an extension
$$0\to \calI_x(1)\to  \calF \to  \calI_y(1)\to 0$$
where $x \ne y$. Since  $\ome = \tome^{**}\not\cong \calO_{\bbP^2}(1)^2$, we conclude that  that $\tome$ is given by a unique non-trivial extension
$$0\to \calI_x(1)\to  \tome \to  \calI_y(1) \to 0,$$
where $x$ is the triple point and $y$ is the intersection point of two lines not passing through $x$. Tensoring by $\calO_L(-2)$ and using that, for any point $z$ , we have an exact sequence
\begin{equation}
0\to \calT or_1^{\bbP^2}(\calO_x,\calO_L) \to \calI_z\otimes_{\calO_{\bbP^2}} \calO_L \to \calO_L(-1)\to 0,
\end{equation}
we see that $W(\tome)$ consists of lines through $x$ or $y$, that is, it is the union of two lines in the dual plane.

Finally $\Ar_{2,5}^{\se}(2)$ consists of arrangements with 2 triple points. The dual set of points lies on the union of two lines, three points on each line, one is the intersection point. The sheaf 
$\Omega^1(\calA)$ is $S$-equivalent to the sheaf $\calO_{\bbP^2}(1)^2$ (in fact, it is isomorphic to this sheaf). It has no jumping lines. The sheaf 
$\tilde{\Omega}^1(\calA)$ is $S$-equivalent to the sheaf $\calI_x(1)\oplus \calI_y(1)$, where $x,y$ are the triple points.  As in the previous case we obtain that $\tome$ is given by a unique non-trivial extension 
$$0\to \calI_x(1)\to \tome\to \calI_y(1) \to 0,$$
where $x,y$ are the triple points of $\calA$. The variety $W(\tome)$ is the union of two lines, dual to the points $x,y$. So, we see that all semi-stable arrangements of 5 lines are not Torelli arrangements. Of course they always lie on a conic.

\end{example}

\begin{example} Let  $m= 6$. In the case when $\calA$ is a generic arrangements the vector bundle $\Omega^1(\calA)$ was extensively studied in \cite{DK2}. Here we are interested in non-generic arrangements.  Since $\mu(\tilde{\Omega}^1(\calA)) = 3/2$, all  sem-stable arrangements are stable. Also we have $\dim \Ar_{2,6} = \dim \calS_{2,6} = 12$, so the map
$$\log:\Ar_{2,6}^{s} \to \calS_{2,6} = \calM_{\bbP^2}(2;3,6) \cong \calM_{\bbP^2}(2;-1,4)$$
is a birational morphism which is an isomorphism on the set of Torelli arrangements.  

 Let $\calA\in  \Ar_{2,6}^{s}(1)$. The bundle $\Omega^1(\calA)$ belongs to the 8-dimensional variety $\calM_{\bbP^2}(2;3,5) \cong \calM_{\bbP^2}(2;-1,3)$. The three lines from $\calA$ which do not pass through the unique triple point $x\in \calA$ are the jumping lines of $\Omega^1(\calA)$. It is known that a  vector bundle $\calE$ from $\calM_{\bbP^2}(2;3,5)$ with 3 non-concurrent jumping lines $L_1,L_2,L_3$ is unique up to an automorphism of $\bbP^2$(\cite{Hu}). Its set of jumping lines is the set $\{L_1,L_2,L_3\}$ and it is given by an extension
\begin{equation}\label{hul}
0\to \calO_{\bbP^2}(1) \to \calE \to \calI_Z(2) \to 0,
\end{equation}
where $Z$ is a 0-dimensional reduced closed subscheme of $\bbP^2$ which consists of three points $p_{ij} = L_i\cap L_j$. Twisting by $\calO_{\bbP^2}(-1)$ we see that 
$$h^0(\calE(-1)) = 1.$$
This shows that the extension is determined uniquely by the isomorphism class of $\calE$. The set of non-isomorphic extensions as in \eqref{hul} is naturally isomorphic to $E =\bbP(H^0(\calO_Z))\cong \bbP^2$. The open supspace of $E$ which consists of section non-vanishing at any point of $Z$ corresponds to stable sheaves. They are all  vector bundles. 
The isomorphism class of $\calE$ is uniquely determined by $Z$ and the class of the extension. Since the map $u\circ \log_1:\Ar_{2,6}^{s}(1) \to \calM_{\bbP^2}(2;3,5)$ is $\PGL(3)$-equivariant, we obtain that any vector bundle from $\calM_{\bbP^2}(2;3,5)$ is isomorphic to $\Omega^1(\calA)$ for some arrangement $\calA$ with $\delta(\calA) = 1$. It determines three lines of $\calA$ not passing through the triple point.

Since any coherent sheaf $\calT$ supported at one point $x$ with $h^0(\calT) = 1$ is isomorphic to the sheaf $\calO_x$,  the sheaf 
$\tilde{\Omega}^1(\calA)$ for such an arrangement  is given by an extension \eqref{second'}
$$0\to \tilde{\Omega}^1(\calA) \to \ome \overset{\alpha}{\to} \calO_x \to 0,$$ 
where $x$ is the triple point of $\calA$. The restriction of $\alpha$ to the subsheaf $\calO_{\bbP^2}(1)$ from \eqref{hul} is not zero. Indeed, otherwise we get that $\tome$ is given by an extension
\begin{equation}\label{ext1}
0\to \calO_{\bbP^2}(1) \to \tome \to \calI_{Z\cup x}(2)\to 0.
\end{equation} 
Tensoring by $\calO_{\bbP^2}(-1)$ we obtain that $h^0(\tome(-1)) = 1$. The residue exact sequence \eqref{res'} shows that $h^0(\tome(-1)) = 0$. In fact, stable sheaves defined  by extensions of type \eqref{ext1} define Hulsbergen vector bundles $\calE$ with $h^0(\calE(-1)) = 1$. They are not isomorphic to $\ome$ for any generic arrangement $\calA$. Since $\alpha$ is not zero on $\calO_{\bbP^2}(1)$ we see that $\tome$ is given by an extension
\begin{equation}\label{ext2}
0\to \calI_x(1) \to \tome \to \calI_Z(2)\to 0,
\end{equation}
where $x$ is the triple point of $\calA$, and $Z$ is the intersection points of the lines not passing through $x$. A standard calculation shows that 
$$\bbP(\Ext_{\bbP^2}^1(\calI_Z(2),\calI_x(1))) \cong \bbP^3.$$ 

Any arrangement of 6 lines with one triple point is a Torelli arrangement. Indeed, suppose $L$ is an unstable line which is not a component of $\calA$.   By tensoring with $\calO_L(-2)$, we easily see that $L$ must contain the triple point. Since $W(\tome)$ cannot be a finite set of more than $6$ point, $W(\tome)$ contains the pencil of lines through $x$. Let $L_1$ be a line from $\calA$ from this pencil. Since the lines $L_2,\ldots,L_6$ form a generic arrangement osculating a nonsingular conic, we have $W(\calA\setminus \{L_1\})$ is the dual conic $C$. By Proposition \ref{new}, $W(\tome) \subset C\cup \{L_1\}$. This shows that $W(\tome)$ cannot contain a line.  Counting parameters we see that any arrangement with one triple point is uniquely determined by the sheaf $\tome$ which is given by a unique extension \eqref{ext2}. So the boundary $\Ar_{2,6}^1$ is birationally isomorphic to a $\bbP^2\times \bbP^1$ fibration over $\calM_{\bbP^2}(2;-1,3)'$, where $\calM_{\bbP^2}(2;-1,3)$ is the open subset of $\calM_{\bbP^2}(2;-1,3)$ representing vector bundles with 3 non-concurrent jumping lines.

Let $\calA\in \Ar_{2,6}^{s}(2)$ be an arrangement with 2 triple points $x,y$. There are two irreducible components of $\Ar_{2,6}^{s}(2)$, each one is  of codimension 2 in $\Ar_{2,6}$. The first one $F_1$ consists of arrangements such that the line $\la x,y\ra$ is a component of $\calA$. The second one $F_2$ consists of arrangements with each line passing through $x$ or $y$. The vector bundle $\ome$ belongs to $\calM_{\bbP^2}(2;3,4) \cong \calM_{\bbP^2}(2;-1,2)$. The variety $\calM_{\bbP^2}(2;-1,2)^0$ is explicitly described in \cite{Hu}. It is isomorphic to the 4-dimensional variety of reducible but not multiple conics. The conic is the conic in $\check{\bbP}^2$ of jumping lines of the second line of a bundle $\calE$ from $\calM_{\bbP^2}(2;3,4)$. Its singular point is the unique jumping line of $\calE$. Each $\calE$ is isomorphic to $\ome$ for some arrangement $\calA$. If $\calA\in F_1$ (resp. $\calA\in F_2$), then the unique jumping line of $\ome$ is the line from $\calA$ which does not pass through the triple points of $\calA$ (resp. the line $\la x,y\ra$)  (see \cite{Sch}). We have an extension
\begin{equation}
0\to \calO_{\bbP^2}(1) \to \ome \to \calI_{Z}(2) \to 0,
\end{equation}
where $Z$ is a closed 0-dimensional subscheme of $\bbP^2$ contained in the jumping line with $h^0(\calO_Z) = 2$. All extension classes with fixed $Z$ are parametrized by $\bbP^1$ and define isomorphic vector bundles. The two points of $Z$ represent the curve jumping lines of the second kind. So, we see that $\ome$ determines very little of $\calA$. 

As in the previous case, one can show that $\tome$ is defined by an extension
\begin{equation}\label{final}
0\to \calI_{x,y}(1) \to \tome \to \calI_{Z}(2) \to 0,
\end{equation}
All such extensions with fixed $Z$ and $x,y$ are parametrized by $\bbP^s$, where 
$s = 3-\#(Z\cap \{x,y\})$. Each isomorphism class of sheaves is determined by a $\bbP^1$ of extensions. 

Any arrangements from $F_1$ is a Torelli arrangement. The proof is similar to the case of arrangements with $\delta(\calA) = 1$. We choose the conic osculating the lines from $\calA$ different from the line $\la x,y\ra$. Thus $F_1$ is a fibration over $\calM_{\bbP^2}(2;-1,2)$ with fibres isomorphic to open subsets of  $\bbP^2\times \bbP^2\times \bbP^2$. The sheaf $\tome$ is given by \eqref{final}, where $Z$ does not lie on the line $\la x,y\ra$.

All arrangements from $F_2$ with the same set of triple points $x,y$ have isomorphic sheaf 
$\tome$. All non-trivial extensions define isomorphic sheaf. By counting constants in parameters defining the shaves in \eqref{final}, we see that the set $Z$ must coincide with the set $\{x,y\}$. The sheaf $\tome$ is a degeneration of the Schwarzenberger bundle corresponding to a nonsingular conic.

Let $\calA\in \Ar_{2,6}^{s}(3)$. The variety $\Ar_{2,6}^{s}(3)$ is an irreducible variety of dimension 8, it belongs to the closure of the irreducible compoinent $F_1$ of $\Ar_{2,6}^{s}(3)$. The arrangement $\calA$ has 3 triple points. In this case 
$\calM_{\bbP^2}(2;3,3)\cong \calM_{\bbP^2}(2;-1,1)$ consists of one point represented by the bundle $\Omega_{\bbP^2}^1(3)$ with no jumping lines. So 
$$\ome \cong  \Omega_{\bbP^2}^1(3) \cong \Theta_{\bbP^2}.$$
 A nonzero section of $\Theta_{\bbP^2}$ defines an extension
 $$0\to \calO_{\bbP^2} \to \Theta_{\bbP^2} \to \calO_{\bbP^2}(3) \to 0.$$
 The sheaf  $\tilde{\Omega}^1(\calA)$ is isomorphic to a kernel of a surjective morphism of sheaves $\ome \to \calO_x\oplus \calO_y\oplus \calO_z$, where $x,y,z$ are the triple points of $\calA$. Arguing as in the previous cases, we obtain that $\tome$  is given by an extension
$$0\to \calI_{x,y,z} \to \tome \to \calO_{\bbP^2}(3) \to 0. $$
The classes of non-trivial extensions are parametrized by $\bbP^2$. The trivial extension is unstable. It is easy to see that any unstable line of $\tome$ must pass through one of the points $x,y,z$, i.e. $W(\tome)$ is contained in the union of three lines. On the other hand, after deleting the line $L = \la x,y\ra$ from $\calA$, we obtain, by Corollary \ref{new} that 
$W(\tome) \subset W(\tilde{\Omega}^1(\calA')\cup \{L\}$, where $\calA'\in \Ar_{2,5}(1)$. It follows from the previous example that the latter consists of two pencils of lines through $z$ and the point $p =L_i\cap L_j$, where $L_i,L_j$ are the lines from $\calA'$ not passing through $z$. Now changing the pair $x,y$ to $x,z$ and $y,z$, and applying the same argument we see that that $\calA$ is a Torelli arrangement.

\end{example} 
Our computations show that the only non-Torelli semi-stable arrangement of 6 lines is the arrangement whose dual points in $\check{\bbP}^2$ are nonsingular points of a conic, nonsingular if the arrangement is generic, and reducible otherwise. This confirms Conjecture \ref{two}.



\begin{thebibliography}{[BPV]}

\bibitem[AO]{AO} V. Ancona, G. Ottaviani, \textit{Unstable hyperplanes for Steiner bundles and multidimensional matrices}, Adv. Geom., {\bf 1} (2001), 165--192.


\bibitem[Ba]{Ba} W. Barth, {\it Moduli of vector bundles on the projective plane}, Invent. Math., {\bf 42} (1977), 63--91.
 
\bibitem[BS]{BS} G. Bohnhorst, H. Spindler, \textit{The stability of certain vector bundles on $\bbP^n$} in ``Complex Algebraic Varieties'', Lect. Notes in Math., {\bf 1507}. Springer-Verlag, 1992, pp. 39--50.

\bibitem[Ca]{Ca} P. Cascini, \textit{On the moduli space of Schwarzenberger bundles}, Pacific J. Math., {\bf 205} (2002), 311--323.

\bibitem[CHKS]{CHKS} F. Catanese, S.  Hosten, A.  Khetan,  B.  Sturmfels, \textit{The maximum likelihood degree}, math.AG/0406533, to appear in Amer. J. Math..



\bibitem[DK1]{DK1}  I. Dolgachev,  M. Kapranov, \textit {Arrangements of hyperplanes and vector bundles on $P^n$}, Duke Math. J., {\bf 71} (1993), 633--664.

\bibitem[DK2]{DK2} I. Dolgachev,  M. Kapranov, \textit{Schur quadrics, cubic surfaces and rank 2 vector bundles over the projective plane}, in ``Journ\'ees de G\'eom\'etrie Alg\'ebrique D'Orsay. Juillet 1992'', Asterisque {\bf 218} (1993), 111--144.

\bibitem[DO]{DO}  I. Dolgachev,  D. Ortland, \textit {Points Sets in Projective Spaces and Theta Functions}, Asterisque., {\bf 165} (1988).


\bibitem[GKZ]{GKZ} I. Gelfand,  M. Kapranov, A. Zelivinsky, \textit{Discriminants, Resultants and Multidimensional Determinants}, Birkh\"auser. 1994.

\bibitem[Ha]{Ha} R. Hartshorne, \textit{Residues and Duality}, Lect. Notes in Math. vol. 20, Springer-Verlag, 1966.

\bibitem[HL]{HL} D. Huybrechts, M. Lehn, \textit{The geometry of moduli space of sheaves}, Aspects of Mathematics, {\bf 31},  1997.

\bibitem[Hu]{Hu} K. Hulek, \textit{Stable rank 2 bundles on $\bbP^2$ with $c_1$ odd}, Math. Ann., {\bf 242}, (1979), 241--266.

\bibitem[Ka]{Ka} S. Karnik, {\it Group actions on moduli of vector bundles}. Ph. D. Thesis, Univ. of Michigan, 2000.
 
 \bibitem[Ma]{Ma} M. Maruyama,  \textit{Moduli of stable sheaves, I and II}, J. Math. Kyoto Univ. {\bf 17} (1977), 91--126, {\bf 18} (1978), 557--614.

\bibitem[MS]{MS} M. Musta\c{t}\u{a}, H. Schenck, \textit{The module of logarithmic p-forms of a locally free arrangement}, J. of Algebra, {\bf 241}, (2001), 699--719.


\bibitem[La]{La} R. Lazarsfeld, \textit{Positivity in Algebraic Geometry II}, Springer. 2004.

\bibitem[Li]{Li} J. Li, {\it Compactification of moduli space of vector bundles over algebraic surfaces}, in ``Collection of papers on geometry, analysis and mathematical physics'', World Sci. Publ., River Edge, 1997, pp. 98--113.
 
\bibitem[OSS]{OSS} C. Okonek, M. Schneider, H. Spindler, \textit{Vector bundles on complex projective space}, Birkh\"auser, 1980.

\bibitem[Ri]{Ri} J.-J. Risler, \textit{Sur l'ideal jacobien d'une courbe plane}, Bull. Soc. Math. France, {\bf 99}, (1971), 305--311.

\bibitem[Sa]{Sa} K.  Saito, \textit{Theory of logarithmic differential forms and logarithmic vector fields}, J. Fac. Sci. Univ. Tokyo, Sci. IA, {\bf 27}, (1980), 265--291.


\bibitem[Sch]{Sch} H. Schenck, \textit{Elementary modifications and line configurations in $\bbP^2$}, Com. Math. Helv, {\bf 78} (2003), 447--462.


\bibitem[Va]{Va} J. Vall\`es, \textit{Nombre maximal d'hyperplanes instables pour un fibr\'e de Steiner}, Math. Zeit. {\bf 233}, 507--514.

\bibitem[Yu]{Yu} S. Yuzvinsky, {\it On generators of the module of logarithmic 1-forms with poles along an arrangement}, J. Alg. Comb., {\bf 4} (1995), 253--269.

\bibitem[Za]{Za} O. Zariski, {\it Characterizations of plane algebroid curves whose module of differentials has maximum torsion}, Proc. Nat. Acad. Sci. USA, {\bf 56} (1966), 781--786.


\bibitem[Zi]{Zi} G. Ziegler, {\it Combinatorial constructions of logarithmic differential forms}, Adv. Math. {\bf 76} (1989), 116--154.




\end{thebibliography}
\end{document}